# LOCAL LINEAR SPATIAL REGRESSION[1]

By Marc Hallin, Zudi Lu and Lanh T. Tran

*Université Libre de Bruxelles, Chinese Academy of Sciences and London School of Economics, and Indiana University*

A local linear kernel estimator of the regression function $\mathbf{x} \mapsto g(\mathbf{x}) := E[Y_\mathbf{i} | \mathbf{X_i} = \mathbf{x}]$, $\mathbf{x} \in \mathbb{R}^d$, of a stationary $(d+1)$-dimensional spatial process $\{(Y_\mathbf{i}, \mathbf{X_i}), \mathbf{i} \in \mathbb{Z}^N\}$ observed over a rectangular domain of the form $\mathcal{I}_\mathbf{n} := \{\mathbf{i} = (i_1, \ldots, i_N) \in \mathbb{Z}^N | 1 \le i_k \le n_k, k = 1, \ldots, N\}$, $\mathbf{n} = (n_1, \ldots, n_N) \in \mathbb{Z}^N$, is proposed and investigated. Under mild regularity assumptions, asymptotic normality of the estimators of $g(\mathbf{x})$ and its derivatives is established. Appropriate choices of the bandwidths are proposed. The spatial process is assumed to satisfy some very general mixing conditions, generalizing classical time-series strong mixing concepts. The size of the rectangular domain $\mathcal{I}_\mathbf{n}$ is allowed to tend to infinity at different rates depending on the direction in $\mathbb{Z}^N$.

**1. Introduction.** Spatial data arise in a variety of fields, including econometrics, epidemiology, environmental science, image analysis, oceanography and many others. The statistical treatment of such data is the subject of an abundant literature, which cannot be reviewed here; for background reading, we refer the reader to the monographs by Anselin and Florax (1995), Cressie (1991), Guyon (1995), Possolo (1991) or Ripley (1981).

Let $\mathbb{Z}^N$, $N \ge 1$, denote the integer lattice points in the $N$-dimensional Euclidean space. A point $\mathbf{i} = (i_1, \ldots, i_N)$ in $\mathbb{Z}^N$ will be referred to as a *site*. Spatial data are modeled as finite realizations of vector stochastic processes indexed by $\mathbf{i} \in \mathbb{Z}^N$: *random fields*. In this paper, we will consider strictly stationary $(d+1)$-dimensional random fields, of the form

(1.1) $$\{(Y_\mathbf{i}, \mathbf{X_i}); \mathbf{i} \in \mathbb{Z}^N\},$$

Received November 2002; revised March 2004.
[1]Supported by a PAI contract of the Belgian Federal Government, an Action de Recherche Concertée of the Communauté Française de Belgique, the National Natural Science Foundation of China, the Leverhulme Trust and NSF Grant DMS-94-03718.

where $Y_\mathbf{i}$, with values in $\mathbb{R}$, and $\mathbf{X}_\mathbf{i}$, with values in $\mathbb{R}^d$, are defined over some probability space $(\Omega, \mathcal{F}, P)$.

A crucial problem for a number of applications is the problem of *spatial regression*, where the influence of a vector $\mathbf{X}_\mathbf{i}$ of covariates on some response variable $Y_\mathbf{i}$ is to be studied in a context of complex spatial dependence. More specifically, assuming that $Y_\mathbf{i}$ has finite expectation, the quantity under study in such problems is the *spatial regression function*

$$g : \mathbf{x} \mapsto g(\mathbf{x}) := E[Y_\mathbf{i} | \mathbf{X}_\mathbf{i} = \mathbf{x}].$$

The spatial dependence structure in this context plays the role of a nuisance, and remains unspecified. Although $g$ of course is only defined up to a P-null set of values of $\mathbf{x}$ (being a class of P-a.s. mutually equal functions rather than a function), we will treat it, for the sake of simplicity, as a well-defined real-valued $\mathbf{x}$-measurable function, which has no implication for the probabilistic statements of this paper. In the particular case under which $\mathbf{X}_\mathbf{i}$ itself is measurable with respect to a subset of $Y_\mathbf{j}$'s, with $\mathbf{j}$ ranging over some neighborhood of $\mathbf{i}$, $g$ is called a *spatial* auto*regression function*. Such spatial autoregression models were considered as early as 1954, in the particular case of a linear autoregression function $g$, by Whittle (1954, 1963); see Besag (1974) for further developments in this context.

In this paper, we are concerned with estimating the spatial regression (autoregression) function $g : \mathbf{x} \mapsto g(\mathbf{x})$; contrary to Whittle (1954), we adopt a nonparametric point of view, avoiding any parametric specification of the possibly extremely complex spatial dependent structure of the data.

For $N = 1$, this problem reduces to the classical problem of (auto)regression for serially dependent observations, which has received extensive attention in the literature; see, for instance, Roussas (1969, 1988), Masry (1983, 1986), Robinson (1983, 1987), Ioannides and Roussas (1987), Masry and Györfi (1987), Yakowitz (1987), Boente and Fraiman (1988), Bosq (1989), Györfi, Härdle, Sarda and Vieu (1989), Tran (1989), Masry and Tjøstheim (1995), Hallin and Tran (1996), Lu and Cheng (1997), Lu (2001) and Wu and Mielniczuk (2002), to quote only a few. Quite surprisingly, despite its importance for applications, the spatial version ($N > 1$) of the same problem remains essentially unexplored. Several recent papers [e.g., Tran (1990), Tran and Yakowitz (1993), Carbon, Hallin and Tran (1996), Hallin, Lu and Tran (2001, 2004), Biau (2003) and Biau and Cadre (2004)] deal with the related problem of estimating the density $f$ of a random field of the form $\{\mathbf{X}_\mathbf{i}; \mathbf{i} \in \mathbb{Z}^N\}$, or the prediction problem but, to the best of our knowledge, the only results available on the estimation of spatial regression functions are those by Lu and Chen (2002, 2004), who investigate the properties of a Nadaraya–Watson kernel estimator for $g$.

Though the Nadaraya–Watson method is central in most nonparametric regression methods in the traditional serial case ($N = 1$), it has been well documented [see, e.g., Fan and Gijbels (1996)] that this approach suffers from several severe drawbacks, such as poor boundary performance, excessive bias and low efficiency, and that the local polynomial fitting methods developed by Stone (1977) and Cleveland (1979) are generally preferable. Local polynomial fitting, and particularly its special case—local linear fitting—recently have become increasingly popular in light of recent work by Cleveland and Loader (1996), Fan (1992), Fan and Gijbels (1992, 1995), Hastie and Loader (1993), Ruppert and Wand (1994) and several others. For $N = 1$, Masry and Fan (1997) have studied the asymptotics of local polynomial fitting for regression under general mixing conditions. In this paper, we extend this approach to the context of spatial regression ($N > 1$) by defining an estimator of $g$ based on local linear fitting and establishing its asymptotic properties.

Extending classical or time-series asymptotics ($N = 1$) to spatial asymptotics ($N > 1$), however, is far from trivial. Due to the absence of any canonical ordering in the space, there is no obvious definition of tail sigma-fields. As a consequence, such a basic concept as ergodicity is all but well defined in the spatial context. And, little seems to exist about this in the literature, where only central limit results are well documented; see, for instance, Bolthausen (1982) or Nakhapetyan (1980). Even the simple idea of a sample size going to infinity (the sample size here is a rectangular domain of the form $\mathcal{I}_\mathbf{n} := \{\mathbf{i} = (i_1, \ldots, i_N) \in \mathbb{Z}^N | 1 \leq i_k \leq n_k, k = 1, \ldots, N\}$, for $\mathbf{n} = (n_1, \ldots, n_N) \in \mathbb{Z}^N$ with strictly positive coordinates $n_1, \ldots, n_N$) or the concept of spatial mixing have to be clarified in this setting. The assumptions we are making (A4), (A4$'$) and (A4$''$) are an attempt to provide reasonable and flexible generalizations of traditional time-series concepts.

Assuming that $\mathbf{x} \mapsto g(\mathbf{x})$ is differentiable at $\mathbf{x}$, with gradient $\mathbf{x} \mapsto g'(\mathbf{x})$, the main idea in local linear regression consists in approximating $g$ in the neighborhood of $\mathbf{x}$ as

$$g(\mathbf{z}) \approx g(\mathbf{x}) + (g'(\mathbf{x}))^\tau (\mathbf{z} - \mathbf{x}),$$

and estimating $(g(\mathbf{x}), g'(\mathbf{x}))$ instead of simply running a classical nonparametric (e.g., kernel-based) estimation method for $g$ itself. In order to do this, we propose a weighted least square estimator $(g_\mathbf{n}(\mathbf{x}), g'_\mathbf{n}(\mathbf{x}))$, and study its asymptotic properties. Mainly, we establish its asymptotic normality under various mixing conditions, as $\mathbf{n}$ goes to infinity in two distinct ways. Either *isotropic divergence* ($\mathbf{n} \Rightarrow \infty$) can be considered; under this case, observations are made over a rectangular domain $\mathcal{I}_\mathbf{n}$ of $\mathbb{Z}^N$ which expands at the same rate in all directions—see Theorems 3.1, 3.2 and 3.5. Or, due to the specific nature of the practical problem under study, the rates of expansion of $\mathcal{I}_\mathbf{n}$ cannot be the same along all directions, and only a less restrictive assumption of possibly *nonisotropic divergence* ($\mathbf{n} \to \infty$) can be made—see Theorems 3.3 and 3.4.

The paper is organized as follows. In Section 2.1 we provide the notation and main assumptions. Section 2.2 introduces the main ideas underlying local linear regression in the context of random fields and sketches the main steps of the proofs to be developed in the sequel. Section 2.3 is devoted to some preliminary results. Section 3 is the main section of the paper, where asymptotic normality is proved under the various types of asymptotics and various mixing assumptions. Section 4 provides some numerical illustrations. Proofs and technical lemmas are concentrated in Section 5.

## 2. Local linear estimation of spatial regression.

2.1. *Notation and main assumptions.* For the sake of convenience, we summarize here the main assumptions we are making on the random field (1.1) and the kernel $K$ to be used in the estimation method. Assumptions (A1)–(A4) are related to the random field itself.

(A1) The random field (1.1) is strictly stationary. For all distinct $\mathbf{i}$ and $\mathbf{j}$ in $\mathbb{Z}^N$, the vectors $\mathbf{X_i}$ and $\mathbf{X_j}$ admit a joint density $f_{\mathbf{i},\mathbf{j}}$; moreover, $|f_{\mathbf{i},\mathbf{j}}(\mathbf{x}', \mathbf{x}'') - f(\mathbf{x}')f(\mathbf{x}'')| \leq C$ for all $\mathbf{i}, \mathbf{j} \in \mathbb{Z}^N$, all $\mathbf{x}', \mathbf{x}'' \in \mathbb{R}^d$, where $C > 0$ is some constant, and $f$ denotes the marginal density of $\mathbf{X_i}$.
(A2) The random variable $Y_{\mathbf{i}}$ has finite absolute moment of order $(2 + \delta)$; that is, $\mathrm{E}[|Y_{\mathbf{i}}|^{2+\delta}] < \infty$ for some $\delta > 0$.
(A3) The spatial regression function $g$ is twice differentiable. Denoting by $g'(\mathbf{x})$ and $g''(\mathbf{x})$ its gradient and the matrix of its second derivatives (at $\mathbf{x}$), respectively, $\mathbf{x} \mapsto g''(\mathbf{x})$ is continuous at all $\mathbf{x}$.

Assumption (A1) is standard in this context; it has been used, for instance, by Masry (1986) in the serial case $N = 1$, and by Tran (1990) in the spatial context $(N > 1)$. If the random field $X_{\mathbf{i}}$ consists of independent observations, then $|f_{\mathbf{i},\mathbf{j}}(\mathbf{x}, \mathbf{x}'') - f(\mathbf{x}')f(\mathbf{x}'')|$ vanishes as soon as $\mathbf{i}$ and $\mathbf{j}$ are distinct. Thus (A1) also allows for unbounded densities.

Assumption (A4) is an assumption of spatial mixing taking two distinct forms [either (A4) and (A4$'$) or (A4) and (A4$''$)]. For any collection of sites $\mathcal{S} \subset \mathbb{Z}^N$, denote by $\mathcal{B}(\mathcal{S})$ the Borel $\sigma$-field generated by $\{(Y_{\mathbf{i}}, \mathbf{X_i}) | \mathbf{i} \in \mathcal{S}\}$; for each couple $\mathcal{S}', \mathcal{S}''$, let $d(\mathcal{S}', \mathcal{S}'') := \min\{\|\mathbf{i}' - \mathbf{i}''\| \,|\, \mathbf{i}' \in \mathcal{S}', \mathbf{i}'' \in \mathcal{S}''\}$ be the distance between $\mathcal{S}'$ and $\mathcal{S}''$, where $\|\mathbf{i}\| := (i_1^2 + \cdots + i_N^2)^{1/2}$ stands for the Euclidean norm. Finally, write $\mathrm{Card}(\mathcal{S})$ for the cardinality of $\mathcal{S}$.

(A4) There exist a function $\varphi$ such that $\varphi(t) \downarrow 0$ as $t \to \infty$, and a function $\psi : \mathbb{N}^2 \to \mathbb{R}^+$ symmetric and decreasing in each of its two arguments, such that the random field (1.1) is mixing, with spatial mixing coefficients $\alpha$ satisfying

$$\alpha(\mathcal{B}(\mathcal{S}'), \mathcal{B}(\mathcal{S}'')) := \sup\{|\mathrm{P}(AB) - \mathrm{P}(A)\mathrm{P}(B)|, A \in \mathcal{B}(\mathcal{S}'), B \in \mathcal{B}(\mathcal{S}'')\}$$
$$\leq \psi(\mathrm{Card}(\mathcal{S}'), \mathrm{Card}(\mathcal{S}''))\varphi(d(\mathcal{S}', \mathcal{S}'')),$$

(2.1)

for any $\mathcal{S}', \mathcal{S}'' \subset \mathbb{Z}^N$. The function $\varphi$, moreover, is such that

$$\lim_{m \to \infty} m^a \sum_{j=m}^{\infty} j^{N-1} \{\varphi(j)\}^{\delta/(2+\delta)} = 0$$

for some constant $a > (4 + \delta)N/(2 + \delta)$.

The assumptions we are making on the function $\psi$ are either

(A4$'$) $\psi(n', n'') \leq \min(n', n'')$

or

(A4$''$) $\psi(n', n'') \leq C(n' + n'' + 1)^\kappa$ for some $C > 0$ and $\kappa > 1$.

In case (2.1) holds with $\psi \equiv 1$, the random field $\{(Y_\mathbf{i}, \mathbf{X}_\mathbf{i})\}$ is called *strongly mixing*.

In the serial case ($N = 1$), many stochastic processes and time series are known to be strongly mixing. Withers (1981) has obtained various conditions for linear processes to be strongly mixing. Under certain weak assumptions, autoregressive and more general nonlinear time-series models are strongly mixing with exponential mixing rates; see Pham and Tran (1985), Pham (1986), Tjøstheim (1990) and Lu (1998). Guyon (1987) has shown that the results of Withers under certain conditions extend to linear random fields, of the form $X_\mathbf{n} = \sum_{\mathbf{j} \in \mathbb{Z}^N} g_\mathbf{j} Z_{\mathbf{n}-\mathbf{j}}$, where the $Z_\mathbf{j}$'s are independent random variables. Assumptions (A4$'$) and (A4$''$) are the same as the mixing conditions used by Neaderhouser (1980) and Takahata (1983), respectively, and are weaker than the uniform strong mixing condition considered by Nakhapetyan (1980). They are satisfied by many spatial models, as shown by Neaderhouser (1980), Rosenblatt (1985) and Guyon (1987).

Throughout, we assume that the random field (1.1) is observed over a rectangular region of the form $\mathcal{I}_\mathbf{n} := \{\mathbf{i} = (i_1, \ldots, i_N) \in \mathbb{Z}^N \mid 1 \leq i_k \leq n_k, k = 1, \ldots, N\}$, for $\mathbf{n} = (n_1, \ldots, n_N) \in \mathbb{Z}^N$ with strictly positive coordinates $n_1, \ldots, n_N$. The total sample size is thus $\hat{\mathbf{n}} := \prod_{k=1}^N n_k$. We write $\mathbf{n} \to \infty$ as soon as $\min_{1 \leq k \leq N} \{n_k\} \to \infty$. The rate at which the rectangular region expands thus can depend on the direction in $\mathbb{Z}^N$. In some problems, however, the assumption that this rate is the same in all directions is natural: we use the notation $\mathbf{n} \Rightarrow \infty$ if $\mathbf{n} \to \infty$ and moreover $|n_j/n_k| < C$ for some $0 < C < \infty$, $1 \leq j, k \leq N$. In this latter case, $\mathbf{n}$ tends to infinity in an *isotropic* way. The *nonisotropic* case $\mathbf{n} \to \infty$ is less restrictive. For more information on the nonisotropic case, we refer to Bradley and Tran (1999) and Lu and Chen (2002).

Assumption (A5) deals with the kernel function $K : \mathbb{R}^d \to \mathbb{R}$ to be used in the estimation method. For any $\mathbf{c} := (c_0, \mathbf{c}_1^\tau)^\tau \in \mathbb{R}^{d+1}$, define

(2.2) $$K_\mathbf{c}(\mathbf{u}) := (c_0 + \mathbf{c}_1^\tau \mathbf{u}) K(\mathbf{u}).$$

(A5)(i) For any $\mathbf{c} \in \mathbb{R}^{d+1}$, $|K_{\mathbf{c}}(\mathbf{u})|$ is uniformly bounded by some constant $K_{\mathbf{c}}^{+}$, and is integrable: $\int_{\mathbb{R}^{d+1}} |K_{\mathbf{c}}(\mathbf{x})| \, d\mathbf{x} < \infty$.

(ii) For any $\mathbf{c} \in \mathbb{R}^{d+1}$, $|K_{\mathbf{c}}|$ has an integrable second-order radial majorant, that is, $Q_{\mathbf{c}}^{K}(\mathbf{x}) := \sup_{\|\mathbf{y}\| \geq \|\mathbf{x}\|}[\|\mathbf{y}\|^{2} K_{\mathbf{c}}(\mathbf{y})]$ is integrable.

Finally, for convenient reference, we list here some conditions on the asymptotic behavior, as $\mathbf{n} \to \infty$, of the bandwidth $b_{\mathbf{n}}$ that will be used in the sequel.

(B1) The bandwidth $b_{\mathbf{n}}$ tends to zero in such a way that $\hat{\mathbf{n}} b_{\mathbf{n}}^{d} \to \infty$ as $\mathbf{n} \to \infty$.

(B2) There exist two sequences of positive integer vectors, $\mathbf{p} = \mathbf{p}_{\mathbf{n}} := (p_1, \ldots, p_N) \in \mathbb{Z}^N$ and $\mathbf{q} = \mathbf{q}_{\mathbf{n}} := (q, \ldots, q) \in \mathbb{Z}^N$, with $q = q_{\mathbf{n}} \to \infty$ such that $p = p_{\mathbf{n}} := \hat{\mathbf{p}} = o((\hat{\mathbf{n}} b_{\mathbf{n}}^{d})^{1/2})$, $q/p_k \to 0$ and $n_k/p_k \to \infty$ for all $k = 1, \ldots, N$, and $\hat{\mathbf{n}} \varphi(q) \to 0$.

(B2′) Same as (B2), but the last condition is replaced by $(\hat{\mathbf{n}}^{\kappa+1}/p)\varphi(q) \to 0$, where $\kappa$ is the constant appearing in (A4″).

(B3) $b_{\mathbf{n}}$ tends to zero in such a manner that $q b_{\mathbf{n}}^{\delta d/[a(2+\delta)]} > 1$ and

$$(2.3) \qquad b_{\mathbf{n}}^{-\delta d/(2+\delta)} \sum_{t=q}^{\infty} t^{N-1} \{\varphi(t)\}^{\delta/(2+\delta)} \to 0 \qquad \text{as } \mathbf{n} \to \infty.$$

2.2. *Local linear fitting.* Local linear fitting consists in approximating, in a neighborhood of $\mathbf{x}$, the unknown function $g$ by a linear function. Under (A3), we have

$$g(\mathbf{z}) \approx g(\mathbf{x}) + (\mathbf{g}'(\mathbf{x}))^{\tau}(\mathbf{z} - \mathbf{x}) := a_0 + \mathbf{a}_1^{\tau}(\mathbf{z} - \mathbf{x}).$$

Locally, this suggests estimating $(a_0, \mathbf{a}_1^{\tau}) = (g(\mathbf{x}), g'(\mathbf{x}))$, hence constructing an estimator of $g$ from

$$(2.4) \quad \begin{pmatrix} g_{\mathbf{n}}(\mathbf{x}) \\ g'_{\mathbf{n}}(\mathbf{x}) \end{pmatrix} = \begin{pmatrix} \hat{a}_0 \\ \hat{\mathbf{a}}_1 \end{pmatrix}$$
$$:= \arg \min_{(a_0, \mathbf{a}_1) \in \mathbb{R}^{d+1}} \sum_{\mathbf{j} \in \mathcal{I}_{\mathbf{n}}} (Y_{\mathbf{j}} - a_0 - \mathbf{a}_1^{\tau}(\mathbf{X}_{\mathbf{j}} - \mathbf{x}))^2 K\left(\frac{\mathbf{X}_{\mathbf{j}} - \mathbf{x}}{b_{\mathbf{n}}}\right),$$

where $b_{\mathbf{n}}$ is a sequence of bandwidths tending to zero at an appropriate rate as $\mathbf{n}$ tends to infinity, and $K(\cdot)$ is a (bounded) kernel with values in $\mathbb{R}^{+}$.

In the classical serial case ($N = 1$; we write $i$ and $n$ instead of $\mathbf{i}$ and $\mathbf{n}$), the solution of the minimization problem (2.4) is easily shown to be $(\mathbf{X}^{\tau} \mathbf{W} \mathbf{X})^{-1} \mathbf{X}^{\tau} \mathbf{W} \mathbf{Y}$, where $\mathbf{X}$ is an $n \times (d+1)$ matrix with $i$th row $(1, b_n^{-1}(\mathbf{X}_i - \mathbf{x})^{\tau})$, $\mathbf{W} = b_n^{-1} \mathrm{diag}(K(\frac{\mathbf{X}_1 - \mathbf{x}}{b_n}), \ldots, K(\frac{\mathbf{X}_n - \mathbf{x}}{b_n}))$, and $\mathbf{Y} = (Y_1, \ldots, Y_n)^{\tau}$ [see, e.g., Fan and Gijbels (1996)]. In the spatial case, things are not as simple, and we rather write the solution to (2.4) as

$$\begin{pmatrix} \hat{a}_0 \\ \hat{\mathbf{a}}_1 b_{\mathbf{n}} \end{pmatrix} = U_{\mathbf{n}}^{-1} V_{\mathbf{n}} \qquad \text{where } \mathbf{V}_{\mathbf{n}} := \begin{pmatrix} v_{\mathbf{n}0} \\ \mathbf{v}_{\mathbf{n}1} \end{pmatrix} \quad \text{and} \quad \mathbf{U}_{\mathbf{n}} := \begin{pmatrix} u_{\mathbf{n}00} & \mathbf{u}_{\mathbf{n}01} \\ \mathbf{u}_{\mathbf{n}10} & \mathbf{u}_{\mathbf{n}11} \end{pmatrix},$$

with [letting $(\frac{\mathbf{X_j}-\mathbf{x}}{b_\mathbf{n}})_0 := 1$]

$$(\mathbf{V_n})_i := (\hat{\mathbf{n}}b_\mathbf{n}^d)^{-1} \sum_{\mathbf{j} \in \mathcal{I}_\mathbf{n}} Y_\mathbf{j} \left(\frac{\mathbf{X_j}-\mathbf{x}}{b_\mathbf{n}}\right)_i K\left(\frac{\mathbf{X_j}-\mathbf{x}}{b_\mathbf{n}}\right), \qquad i = 0, \ldots, d,$$

and

$$(\mathbf{U_n})_{i\ell} := (\hat{\mathbf{n}}b_\mathbf{n}^d)^{-1} \sum_{\mathbf{j} \in \mathcal{I}_\mathbf{n}} \left(\frac{\mathbf{X_j}-\mathbf{x}}{b_\mathbf{n}}\right)_i \left(\frac{\mathbf{X_j}-\mathbf{x}}{b_\mathbf{n}}\right)_\ell K\left(\frac{\mathbf{X_j}-\mathbf{x}}{b_\mathbf{n}}\right), \qquad i, \ell = 0, \ldots, d.$$

It follows that

(2.5)
$$\begin{aligned}\mathbf{H_n} &:= \begin{pmatrix} \hat{a}_0 - a_0 \\ \hat{\mathbf{a}}_1 b_\mathbf{n} - \mathbf{a}_1 b_\mathbf{n} \end{pmatrix} = \begin{pmatrix} g_\mathbf{n}(\mathbf{x}) - g(\mathbf{x}) \\ (g'_\mathbf{n}(\mathbf{x}) - g'(\mathbf{x}))b_\mathbf{n} \end{pmatrix} \\ &= \mathbf{U}_\mathbf{n}^{-1}\left\{\mathbf{V_n} - \mathbf{U_n}\begin{pmatrix} a_0 \\ \mathbf{a}_1 b_\mathbf{n} \end{pmatrix}\right\} =: \mathbf{U}_\mathbf{n}^{-1}\mathbf{W_n},\end{aligned}$$

where

(2.6)
$$\mathbf{W_n} := \begin{pmatrix} w_{\mathbf{n}0} \\ \mathbf{w}_{\mathbf{n}1} \end{pmatrix},$$
$$(\mathbf{W_n})_i := (\hat{\mathbf{n}}b_\mathbf{n}^d)^{-1} \sum_{\mathbf{j} \in \mathcal{I}_\mathbf{n}} Z_\mathbf{j}\left(\frac{\mathbf{X_j}-\mathbf{x}}{b_\mathbf{n}}\right)_i K\left(\frac{\mathbf{X_j}-\mathbf{x}}{b_\mathbf{n}}\right), \qquad i = 0, \ldots, d,$$

and $Z_\mathbf{j} := Y_\mathbf{j} - a_0 - \mathbf{a}_1^\tau(\mathbf{X_j}-\mathbf{x})$.

The organization of the paper is as follows. If, under adequate conditions, we are able to show that:

(C1) $(\hat{\mathbf{n}}b_\mathbf{n}^d)^{1/2}(\mathbf{W_n} - E\mathbf{W_n})$ is asymptotically normal,
(C2) $(\hat{\mathbf{n}}b_\mathbf{n}^d)^{1/2}E\mathbf{W_n} \to \mathbf{0}$ and $\text{Var}((\hat{\mathbf{n}}b_\mathbf{n}^d)^{1/2}\mathbf{W_n}) \to \boldsymbol{\Sigma}$, and
(C3) $\mathbf{U_n} \xrightarrow{P} \mathbf{U}$,

then (2.5) and Slutsky's classical argument imply that, for all $\mathbf{x}$ (all quantities involved indeed depend on $\mathbf{x}$),

$$(\hat{\mathbf{n}}b_\mathbf{n}^d)^{1/2}\begin{pmatrix} g_\mathbf{n}(\mathbf{x}) - g(\mathbf{x}) \\ (g'_\mathbf{n}(\mathbf{x}) - g'(\mathbf{x}))b_\mathbf{n} \end{pmatrix} = (\hat{\mathbf{n}}b_\mathbf{n}^d)^{1/2}\mathbf{H_n} \xrightarrow{\mathcal{L}} \mathcal{N}(\mathbf{0}, \mathbf{U}^{-1}\boldsymbol{\Sigma}(\mathbf{U}^{-1})^\tau).$$

This asymptotic normality result (with explicit values of $\boldsymbol{\Sigma}$ and $\mathbf{U}$), under various forms (depending on the mixing assumptions [(A4$'$) or (A4$''$)], the choice of the bandwidth $b_\mathbf{n}$, the way $\mathbf{n}$ tends to infinity, etc.), is the main contribution of this paper; see Theorems 3.1–3.5. Section 2.3 deals with (C2) and (C3) under $\mathbf{n} \to \infty$ (hence also under the stronger assumption that $\mathbf{n} \Rightarrow \infty$), and Sections 3.1 and 3.2 with (C1) under $\mathbf{n} \Rightarrow \infty$ and $\mathbf{n} \to \infty$, respectively.

2.3. *Preliminaries.* Claim (C3) is easily established from the following lemma, the proof of which is similar to that of Lemma 2.2, and is therefore omitted.

LEMMA 2.1. *Assume that* (A1), (A4) *and* (A5) *hold, that* $b_\mathbf{n}$ *satisfies assumption* (B1) *and that* $n_k b_\mathbf{n}^{\delta d/[a(2+\delta)]} > 1$ *as* $\mathbf{n} \to \infty$. *Then, for all* $\mathbf{x}$,

$$\mathbf{U_n} \xrightarrow{P} \mathbf{U} := \begin{pmatrix} f(\mathbf{x}) \int K(\mathbf{u})\,d\mathbf{u} & f(\mathbf{x}) \int \mathbf{u}^\tau K(\mathbf{u})\,d\mathbf{u} \\ f(\mathbf{x}) \int \mathbf{u} K(\mathbf{u})\,d\mathbf{u} & f(\mathbf{x}) \int \mathbf{u}\mathbf{u}^\tau K(\mathbf{u})\,d\mathbf{u} \end{pmatrix}$$

*as* $\mathbf{n} \to \infty$.

The remainder of this section is devoted to claim (C2). The usual Cramér–Wold device will be adopted. For all $\mathbf{c} := (c_0, \mathbf{c}_1^\tau)^\tau \in \mathbb{R}^{1+d}$, let

$$A_\mathbf{n} := (\hat{\mathbf{n}} b_\mathbf{n}^d)^{1/2} \mathbf{c}^\tau \mathbf{W_n} = (\hat{\mathbf{n}} b_\mathbf{n}^d)^{-1/2} \sum_{\mathbf{j} \in \mathcal{I}_\mathbf{n}} Z_\mathbf{j} K_\mathbf{c}\left(\frac{\mathbf{X_j} - \mathbf{x}}{b_\mathbf{n}}\right),$$

with $K_\mathbf{c}(\mathbf{u})$ defined in (2.2). The following lemma provides the asymptotic variance of $A_\mathbf{n}$ for all $\mathbf{c}$, hence that of $(\hat{\mathbf{n}} b_\mathbf{n}^d)^{1/2} \mathbf{W_n}$.

LEMMA 2.2. *Assume that* (A1), (A2), (A4) *and* (A5) *hold, that* $b_\mathbf{n}$ *satisfies assumption* (B1) *and that* $n_k b_\mathbf{n}^{\delta d/[(2+\delta)a]} > 1$ *for all* $k = 1,\ldots,N$, *as* $\mathbf{n} \to \infty$. *Then*

$$(2.7) \quad \lim_{\mathbf{n} \to \infty} \mathrm{Var}[A_\mathbf{n}] = \mathrm{Var}(Y_\mathbf{j}|\mathbf{X_j} = \mathbf{x}) f(\mathbf{x}) \int_{\mathbb{R}^d} K_\mathbf{c}^2(\mathbf{u})\,d\mathbf{u} = \mathbf{c}^\tau \mathbf{\Sigma} \mathbf{c},$$

*where*

$$\mathbf{\Sigma} := \mathrm{Var}(Y_\mathbf{j}|\mathbf{X_j} = \mathbf{x}) f(\mathbf{x}) \begin{pmatrix} \int K^2(\mathbf{u})\,d\mathbf{u} & \int \mathbf{u}^\tau K^2(\mathbf{u})\,d\mathbf{u} \\ \int \mathbf{u} K^2(\mathbf{u})\,d\mathbf{u} & \int \mathbf{u}\mathbf{u}^\tau K^2(\mathbf{u})\,d\mathbf{u} \end{pmatrix}.$$

*Hence* $\lim_{\mathbf{n} \to \infty} \mathrm{Var}((\hat{\mathbf{n}} b_\mathbf{n}^d)^{1/2} \mathbf{W_n}) = \mathbf{\Sigma}$.

For the proof see Section 5.1.

Next we consider the asymptotic behavior of $\mathrm{E}[A_\mathbf{n}]$.

LEMMA 2.3. *Under assumptions* (A3) *and* (A5),

$$(2.8) \quad \begin{aligned} \mathrm{E}[A_\mathbf{n}] &= \sqrt{\hat{\mathbf{n}} b_\mathbf{n}^d}\, b_\mathbf{n}^2 \tfrac{1}{2} f(\mathbf{x})\, \mathrm{tr}\left[g''(\mathbf{x}) \int \mathbf{u}\mathbf{u}^\tau K_\mathbf{c}(\mathbf{u})\,d\mathbf{u}\right] + o(\sqrt{\hat{\mathbf{n}} b_\mathbf{n}^d}\, b_\mathbf{n}^2) \\ &= \sqrt{\hat{\mathbf{n}} b_\mathbf{n}^d}\, b_\mathbf{n}^2 [c_0 B_0(\mathbf{x}) + \mathbf{c}_1^\tau \mathbf{B}_1(\mathbf{x})] + o(\sqrt{\hat{\mathbf{n}} b_\mathbf{n}^d}\, b_\mathbf{n}^2), \end{aligned}$$

*where*

$$B_0(\mathbf{x}) := \tfrac{1}{2} f(\mathbf{x}) \sum_{i=1}^{d} \sum_{j=1}^{d} g_{ij}(\mathbf{x}) \int u_i u_j K(\mathbf{u}) \, d\mathbf{u},$$

$$\mathbf{B}_1(\mathbf{x}) := \tfrac{1}{2} f(\mathbf{x}) \sum_{i=1}^{d} \sum_{j=1}^{d} g_{ij}(\mathbf{x}) \int u_i u_j \mathbf{u} K(\mathbf{u}) \, d\mathbf{u},$$

$g_{ij}(\mathbf{x}) = \partial^2 g(\mathbf{x})/\partial x_i \, \partial x_j$, $i, j = 1, \ldots, d$, *and* $\mathbf{u} := (u_1, \ldots u_d)^\tau \in \mathbb{R}^d$.

For the proof see Section 5.2.

### 3. Asymptotic normality.

3.1. *Asymptotic normality under mixing assumption* (A4′). The asymptotic normality of our estimators relies in a crucial manner on the following lemma [see (2.6) for the definition of $\mathbf{W_n}(\mathbf{x})$].

LEMMA 3.1. *Suppose that assumptions* (A1), (A2), (A4), (A4′) *and* (A5) *hold, and that the bandwidth* $b_\mathbf{n}$ *satisfies conditions* (B1)–(B3). *Denote by* $\sigma^2$ *the asymptotic variance* (2.7). *Then* $(\hat{\mathbf{n}} b_\mathbf{n}^d)^{1/2} (\mathbf{c}^\tau [\mathbf{W_n}(\mathbf{x}) - \mathrm{E}\mathbf{W_n}(\mathbf{x})]/\sigma)$ *is asymptotically standard normal as* $\mathbf{n} \to \infty$.

For the proof see Section 5.3.

We now turn to the main consistency and asymptotic normality results. First, we consider the case where the sample size tends to $\infty$ in the manner of Tran (1990), that is, $\mathbf{n} \Rightarrow \infty$.

THEOREM 3.1. *Let assumptions* (A1)–(A3), (A4′) *and* (A5) *hold, with* $\varphi(x) = O(x^{-\mu})$ *for some* $\mu > 2(3+\delta)N/\delta$. *Suppose that there exists a sequence of positive integers* $q = q_\mathbf{n} \to \infty$ *such that* $q_\mathbf{n} = o((\hat{\mathbf{n}} b_\mathbf{n}^d)^{1/(2N)})$ *and* $\hat{\mathbf{n}} q^{-\mu} \to 0$ *as* $\mathbf{n} \Rightarrow \infty$, *and that the bandwidth* $b_\mathbf{n}$ *tends to zero in such a manner that*

(3.1) $$q b_\mathbf{n}^{\delta d/[a(2+\delta)]} > 1$$

*for some* $(4+\delta)N/(2+\delta) < a < \mu\delta/(2+\delta) - N$ *as* $\mathbf{n} \Rightarrow \infty$. *Then*,

(3.2) $$(\hat{\mathbf{n}} b_\mathbf{n}^d)^{1/2} \left[ \begin{pmatrix} g_\mathbf{n}(\mathbf{x}) - g(\mathbf{x}) \\ b_\mathbf{n}(g'_\mathbf{n}(\mathbf{x}) - g'(\mathbf{x})) \end{pmatrix} - \mathbf{U}^{-1} \begin{pmatrix} B_0(\mathbf{x}) \\ \mathbf{B}_1(\mathbf{x}) \end{pmatrix} b_\mathbf{n}^2 \right] \xrightarrow{\mathcal{L}} \mathcal{N}(\mathbf{0}, \mathbf{U}^{-1} \mathbf{\Sigma} (\mathbf{U}^{-1})^\tau)$$

*as* $\mathbf{n} \Rightarrow \infty$, *where* $\mathbf{U}$, $\mathbf{\Sigma}$, $B_0(\mathbf{x})$ *and* $\mathbf{B}_1(\mathbf{x})$ *are defined in Lemmas* 2.1, 2.2 *and* 2.3, *respectively. If, furthermore, the kernel* $K(\cdot)$ *is a symmetric density function, then* (3.2) *can be reinforced into*

$$\begin{pmatrix} (\hat{\mathbf{n}} b_\mathbf{n}^d)^{1/2} [g_\mathbf{n}(\mathbf{x}) - g(\mathbf{x}) - B_g(\mathbf{x}) b_\mathbf{n}^2] \\ (\hat{\mathbf{n}} b_\mathbf{n}^{d+2})^{1/2} [g'_\mathbf{n}(\mathbf{x}) - g'(\mathbf{x})] \end{pmatrix} \xrightarrow{\mathcal{L}} \mathcal{N}\left( \mathbf{0}, \begin{pmatrix} \sigma_0^2(\mathbf{x}) & 0 \\ 0 & \boldsymbol{\sigma}_1^2(\mathbf{x}) \end{pmatrix} \right)$$

[*so that $g_\mathbf{n}(\mathbf{x})$ and $g'_\mathbf{n}(\mathbf{x})$ are asymptotically independent*], where

$$B_g(\mathbf{x}) := \tfrac{1}{2} \sum_{i=1}^d g_{ii}(\mathbf{x}) \int (\mathbf{u})_i^2 K(\mathbf{u}) \, d\mathbf{u}, \qquad \sigma_0^2(\mathbf{x}) := \frac{\operatorname{Var}(Y_\mathbf{j} | \mathbf{X}_\mathbf{j} = \mathbf{x}) \int K^2(\mathbf{u}) \, d\mathbf{u}}{f(\mathbf{x})}$$

and

$$\boldsymbol{\sigma}_1^2(\mathbf{x}) := \frac{\operatorname{Var}(Y_\mathbf{j} | \mathbf{X}_\mathbf{j} = \mathbf{x})}{f(\mathbf{x})}$$
$$\times \left[ \int \mathbf{u}\mathbf{u}^\tau K(\mathbf{u}) \, d\mathbf{u} \right]^{-1} \left[ \int \mathbf{u}\mathbf{u}^\tau K^2(\mathbf{u}) \, d\mathbf{u} \right] \left[ \int \mathbf{u}\mathbf{u}^\tau K(\mathbf{u}) \, d\mathbf{u} \right]^{-1}.$$

The asymptotic normality results in Theorem 3.1 are stated for $g_\mathbf{n}(\mathbf{x})$ and $g'_\mathbf{n}(\mathbf{x})$ at a given site $\mathbf{x}$. They are easily extended, via the traditional Cramér–Wold device, into a joint asymptotic normality result for any couple $(\mathbf{x}_1, \mathbf{x}_2)$ (or any finite collection) of sites; the asymptotic covariance terms [between $g_\mathbf{n}(\mathbf{x}_1)$ and $g_\mathbf{n}(\mathbf{x}_2)$, $g_\mathbf{n}(\mathbf{x}_1)$ and $g'_\mathbf{n}(\mathbf{x}_2)$, etc.] all are equal to zero, as in related results on density estimation [see Hallin and Tran (1996) or Lu (2001)]. The same remark also holds for Theorems 3.2–3.5 below.

PROOF OF THEOREM 3.1. Since $q$ is $o((\hat{\mathbf{n}} b_\mathbf{n}^d)^{1/2N})$, there exists $s_\mathbf{n} \to 0$ such that $q = (\hat{\mathbf{n}} b_\mathbf{n}^d)^{1/2N} s_\mathbf{n}$. Take $p_k := (\hat{\mathbf{n}} b_\mathbf{n}^d)^{1/2N} s_\mathbf{n}^{1/2}$, $k = 1, \ldots, N$. Then $q/p_k = s_\mathbf{n}^{1/2} \to 0$, $\hat{\mathbf{p}} = (\hat{\mathbf{n}} b_\mathbf{n}^d)^{1/2} s_\mathbf{n}^{N/2} = o((\hat{\mathbf{n}} b_\mathbf{n}^d)^{1/2})$ and $\hat{\mathbf{n}} \varphi(q) = \hat{\mathbf{n}} q^{-\mu} \to 0$. As $\mathbf{n} \Rightarrow \infty$, $p := \hat{\mathbf{p}} < (\hat{\mathbf{n}} b_\mathbf{n}^d)^{1/2}$ for large $\hat{\mathbf{n}}$. It follows that $\hat{\mathbf{n}}/p > (\hat{\mathbf{n}} b_\mathbf{n}^{-d})^{1/2} \to \infty$, hence $n_k/p_k \to \infty$ for all $k$. Thus, condition (B2) is satisfied.

Because $\varphi(j) = Cj^{-\mu}$,

$$m^a \sum_{j=m}^\infty j^{N-1} \{\varphi(j)\}^{\delta/(2+\delta)} = Cm^a \sum_{j=m}^\infty j^{N-1} j^{-\mu\delta/(2+\delta)}$$
$$\leq Cm^a m^{N-\mu\delta/(2+\delta)} = m^{-[\mu\delta/(2+\delta)-a-N]},$$

a quantity that tends to zero as $m \to \infty$ since $(4+\delta)N/(2+\delta) < a < \mu\delta/(2+\delta) - N$, hence $\mu\delta/(2+\delta) > a + N$. Assumption (A4) and the fact that $q b_\mathbf{n}^{\delta d/[a(2+\delta)]} > 1$ imply that $b_\mathbf{n}^{-\delta d/(2+\delta)} < q^a$ and that (2.3) holds. Now

$$\mathbf{H_n} - \mathbf{U}^{-1} \mathbf{E} \mathbf{W_n} = \mathbf{U_n}^{-1}(\mathbf{W_n} - \mathbf{E} \mathbf{W_n}) + (\mathbf{U_n}^{-1} - \mathbf{U}^{-1}) \mathbf{E} \mathbf{W_n}.$$

The theorem thus follows from Lemmas 2.1, 2.3 and 3.1. □

One of the important advantages of local polynomial (and linear) fitting over the more traditional Nadaraya–Watson approach is that it has much better boundary behavior. This advantage often has been emphasized in the usual regression and time-series settings when the regressors take values on a compact subset of $\mathbb{R}^d$. For example, as Fan and Gijbels (1996) and Fan and Yao (2003) illustrate, for a univariate regressor $X$ with bounded support ($[0, 1]$, say; here, $d = 1$), it can be proved, using an argument similar to the one we develop in the proof of Theorem 3.1, that asymptotic normality still holds at the boundary point $x = cb_{\mathbf{n}}$ (here $c$ is a positive constant), but with asymptotic bias and variances

$$
\begin{aligned}
B_g &:= \frac{1}{2}\left(\frac{\partial^2 g}{\partial x^2}\right)_{x=0^+} \int_{-c}^{\infty} u^2 K(u)\, du, \\
\sigma_0^2 &:= \frac{\operatorname{Var}(Y_{\mathbf{j}}|X_{\mathbf{j}} = 0^+) \int_{-c}^{\infty} K^2(u)\, du}{f(0^+)}
\end{aligned}
\tag{3.3}
$$

and

$$
\sigma_1^2 := \frac{\operatorname{Var}(Y_{\mathbf{j}}|X_{\mathbf{j}} = 0^+)}{f(0^+)} \left[\int_{-c}^{\infty} u^2 K(u)\, du\right]^{-2} \left[\int_{-c}^{\infty} u^2 K^2(u)\, du\right],
\tag{3.4}
$$

respectively. This advantage is likely to be much more substantial as $N$ grows. Therefore, results on the model of (3.3) and (3.4) on the boundary behavior of our estimators would be highly desirable. Such results, however, are all but straightforward, and we leave them for future research. On the other hand, the statistical relevance of boundary effects is also of lesser importance, as the ultimate objective in random fields, as opposed to time series, seldom consists in "forecasting" the process beyond the boundary of the observed domain.

In the important particular case under which $\varphi(x)$ tends to zero at an exponential rate, the same results are obtained under milder conditions.

THEOREM 3.2. *Let assumptions* (A1)–(A3), (A4′) *and* (A5) *hold, with* $\varphi(x) = O(e^{-\xi x})$ *for some* $\xi > 0$. *Then, if* $b_{\mathbf{n}}$ *tends to zero as* $\mathbf{n} \Rightarrow \infty$ *in such a manner that*

$$
(\hat{\mathbf{n}} b_{\mathbf{n}}^{d(1+2N\delta/a(2+\delta))})^{1/2N}(\log \hat{\mathbf{n}})^{-1} \to \infty
\tag{3.5}
$$

*for some* $a > (4+\delta)N/(2+\delta)$, *the conclusions of Theorem* 3.1 *still hold.*

PROOF. By (3.5), there exists a monotone positive function $\mathbf{n} \mapsto g(\mathbf{n})$ such that $g(\mathbf{n}) \to \infty$ and $(\hat{\mathbf{n}} b_{\mathbf{n}}^{d(1+2N\delta/a(2+\delta))})^{1/2N}(g(\mathbf{n})\log \hat{\mathbf{n}})^{-1} \to \infty$ as $\mathbf{n} \Rightarrow \infty$. Let $q := (\hat{\mathbf{n}} b_{\mathbf{n}}^d)^{1/2N}(g(\mathbf{n}))^{-1}$, and $p_k := (\hat{\mathbf{n}} b_{\mathbf{n}}^d)^{1/2N} g^{-1/2}(\mathbf{n})$. Then $q/p_k = g^{-1/2}(\mathbf{n}) \to 0$, $\hat{\mathbf{p}} = (\hat{\mathbf{n}} b_{\mathbf{n}}^d)^{1/2} g^{-N/2}(\mathbf{n}) = o((\hat{\mathbf{n}} b_{\mathbf{n}}^d)^{1/2})$ and $n_k/p_k \to \infty$ as $\mathbf{n} \Rightarrow \infty$. For arbitrary $C > 0$, $q \geq C \log \hat{\mathbf{n}}$ for sufficiently large $\hat{\mathbf{n}}$. Thus

$$
\hat{\mathbf{n}} \varphi(q) \leq C \hat{\mathbf{n}} e^{-\xi q} \leq C \hat{\mathbf{n}} \exp(-C\xi \log \hat{\mathbf{n}}) = C \hat{\mathbf{n}}^{-C\xi+1},
$$

which tends to zero if we choose $C > 1/\xi$. Hence condition (B2) is satisfied. Next, for $0 < \xi' < \xi$,

$$q^a \sum_{i=q}^{\infty} i^{N-1} \varphi(i)^{\delta/(2+\delta)} \leq Cq^a \sum_{i=q}^{\infty} i^{N-1} e^{-\xi i \delta/(2+\delta)}$$

$$\leq Cq^a \sum_{i=q}^{\infty} e^{-\xi' i \delta/(2+\delta)}$$

$$\leq Cq^a e^{-\xi' q \delta/(2+\delta)}.$$

Note that $b_{\mathbf{n}}^d \geq C\hat{\mathbf{n}}^{-1}$ and $q > C \log \hat{\mathbf{n}}$, so that assumption (A4) holds. In addition,

$$qb_{\mathbf{n}}^{\delta d/[a(2+\delta)]} = (\hat{\mathbf{n}} b_{\mathbf{n}}^{d+2N\delta d/a(2+\delta)})^{1/2N} (g(\mathbf{n}))^{-1} > 1$$

for $\hat{\mathbf{n}}$ large enough. It is easily verified that this implies that condition (B3) is satisfied. The theorem follows. □

Note that, in the one-dimensional case $N = 1$, and for "large" values of $a$, the condition (3.5) is "close" to the condition that $nb_n^d \to \infty$, which is usual in the classical case of independent observations.

Next we consider the situation under which the sample size tends to $\infty$ in the "weak" sense (i.e., $\mathbf{n} \to \infty$ instead of $\mathbf{n} \Rightarrow \infty$).

THEOREM 3.3. *Let assumptions* (A1)–(A3), (A4') *and* (A5) *hold, with* $\varphi(x) = O(x^{-\mu})$ *for some* $\mu > 2(3+\delta)N/\delta$. *Let the sequence of positive integers* $q = q_\mathbf{n} \to \infty$ *and the bandwidth* $b_\mathbf{n}$ *factor into* $b_\mathbf{n} := \prod_{i=1}^{N} b_{n_i}$, *such that* $\hat{\mathbf{n}} q^{-\mu} \to 0$, $q = o(\min_{1 \leq k \leq N}(n_k b_{n_k}^d)^{1/2})$, *and*

$$qb_{\mathbf{n}}^{\delta d/a(2+\delta)} > 1 \qquad \textit{for some } (4+\delta)N/(2+\delta) < a < \mu\delta/(2+\delta) - N.$$

*Then the conclusions of Theorem* 3.1 *hold as* $\mathbf{n} \to \infty$.

PROOF. Since $q = o(\min_{1 \leq k \leq N}(n_k b_{n_k}^d)^{1/2})$, there exists a sequence $s_{n_k} \to 0$ such that

$$q = \min_{1 \leq k \leq N}((n_k b_{n_k}^d)^{1/2} s_{n_k}) \qquad \text{as } \mathbf{n} \to \infty.$$

Take $p_k = (n_k b_{n_k}^d)^{1/2} s_{n_k}^{1/2}$. Then $q/p_k \leq s_{n_k}^{1/2} \to 0$, $\hat{\mathbf{p}} = (\hat{\mathbf{n}} b_\mathbf{n}^d)^{1/2} \prod_{k=1}^{N} s_{n_k}^{1/2} = o((\hat{\mathbf{n}} b_\mathbf{n}^d)^{1/2})$ and $\hat{\mathbf{n}} \varphi(q) = \hat{\mathbf{n}} q^{-\mu} \to 0$. As $\mathbf{n} \to \infty$, $p_k < (n_k b_{n_k}^d)^{1/2}$, hence $n_k/p_k > (n_k b_{n_k}^{-d})^{1/2} \to \infty$. Thus condition (B2) is satisfied. The end of the proof is entirely similar to that of Theorem 3.1. □

In the important case that $\varphi(x)$ tends to zero at an exponential rate, we have the following result, which parallels Theorem 3.2.

THEOREM 3.4. *Let assumptions* (A1)–(A3), (A4′) *and* (A5) *hold, with* $\varphi(x) = O(e^{-\xi x})$ *for some* $\xi > 0$. *Let the bandwidth* $b_{\mathbf{n}}$ *factor into* $b_{\mathbf{n}} := \prod_{i=1}^{N} b_{n_i}$ *in such a way that, as* $\mathbf{n} \to \infty$,

$$\min_{1 \leq k \leq N} \{(n_k b_{n_k}^d)^{1/2}\} b_{\mathbf{n}}^{d\delta/a(2+\delta)} (\log \hat{\mathbf{n}})^{-1} \to \infty \tag{3.6}$$

*for some* $a > (4+\delta)N/(2+\delta)$. *Then the conclusions of Theorem* 3.1 *hold as* $\mathbf{n} \to \infty$.

PROOF. By (3.6) there exist positive sequences indexed by $n_k$ such that $g_{n_k} \uparrow \infty$ as $n_k \to \infty$ and

$$\min_{1 \leq k \leq N} \{(n_k b_{n_k}^d)^{1/2} g_{n_k}^{-1}\} b_{\mathbf{n}}^{d\delta/a(2+\delta)} (\log \hat{\mathbf{n}})^{-1} \to \infty$$

as $\mathbf{n} \to \infty$. Let $q := \min_{1 \leq k \leq N} \{(n_k b_{n_k}^d)^{1/2} (g_{n_k})^{-1}\}$ and $p_k := (n_k b_{n_k}^d)^{1/2} g_{n_k}^{-1/2}$. Then $q/p_k \leq g_{n_k}^{-1/2} \to 0$, $\hat{\mathbf{p}} = (\hat{\mathbf{n}} b_{\mathbf{n}}^d)^{1/2} \prod_{k=1}^{N} g_{n_k}^{-1/2} = o((\hat{\mathbf{n}} b_{\mathbf{n}}^d)^{1/2})$ and $n_k/p_k = (n_k b_{n_k}^{-d})^{1/2} g_{n_k}^{1/2} \to \infty$ as $\mathbf{n} \to \infty$. For arbitrary $C > 0$, $q \geq C \log \hat{\mathbf{n}}$ for sufficiently large $\hat{\mathbf{n}}$. Thus

$$\hat{\mathbf{n}} \varphi(q) \leq C \hat{\mathbf{n}} e^{-\xi q} \leq C \hat{\mathbf{n}} \exp(-C \xi \log \hat{\mathbf{n}}) = C \hat{\mathbf{n}}^{-C\xi + 1},$$

which tends to zero for $C > 1/\xi$. Hence, condition (B2) is satisfied. Next, for $0 < \xi' < \xi$,

$$q^a \sum_{i=q}^{\infty} i^{N-1} \varphi(i)^{\delta/(2+\delta)} \leq C q^a \sum_{i=q}^{\infty} i^{N-1} e^{-\xi i \delta/(2+\delta)}$$

$$\leq C q^a \sum_{i=q}^{\infty} e^{-\xi' i \delta/(2+\delta)}$$

$$\leq C q^a e^{-\xi' q \delta/(2+\delta)}.$$

Note that $q > C \log \hat{\mathbf{n}}$. Assumption (A4′) and (3.1) imply that $q b_{\mathbf{n}}^{\delta d/a(2+\delta)} > 1$ for $\mathbf{n}$ large enough. This in turn implies that condition (B3) is satisfied. The theorem follows. □

3.2. *Asymptotic normality under mixing assumption* (A4″). We start with an equivalent, under (A4″), of Lemma 3.1.

LEMMA 3.2. *Suppose that assumptions* (A1), (A2), (A4) *or* (A4″), *and* (A5) *hold, and that the bandwidth* $b_{\mathbf{n}}$ *satisfies conditions* (B1), (B2′) *and* (B3). *Then the conclusions of Lemma* 3.1 *still hold as* $\mathbf{n} \to \infty$.

PROOF. The proof is a slight variation of the argument of Lemma 3.1, and we describe it only briefly. The only significant difference is in the checking of (5.18). Let $\widetilde{U}_1, \ldots, \widetilde{U}_M$ be as in Lemma 3.1. By Lemma 5.3 and assumption (A4″),

$$Q_1 \leq C \sum_{i=1}^{M} [\hat{\mathbf{p}} + (M-i)\hat{\mathbf{p}} + 1]^\kappa \varphi(q)$$
$$\leq C\hat{\mathbf{p}}^\kappa M^{\kappa+1} \varphi(q) \leq C(\hat{\mathbf{n}}^{(\kappa+1)}/\hat{\mathbf{p}})\varphi(q),$$

which tends to zero by condition (B2′); (5.18) follows. □

We then have the following counterpart of Theorem 3.1.

THEOREM 3.5. *Let assumptions* (A1)–(A3), (A4″) *and* (A5) *hold, with* $\varphi(x) = O(x^{-\mu})$ *for some* $\mu > 2(3+\delta)N/\delta$. *Suppose that there exists a sequence of positive integers* $q = q_\mathbf{n} \to \infty$ *such that* $q_\mathbf{n} = o((\hat{\mathbf{n}} b_\mathbf{n}^d)^{1/2N})$ *and* $\hat{\mathbf{n}}^{\kappa+1} q^{-\mu-N} \to 0$ *as* $\mathbf{n} \Rightarrow \infty$, *and that the bandwidth* $b_\mathbf{n}$ *tends to zero in such a manner that* (3.1) *is satisfied as* $\mathbf{n} \Rightarrow \infty$. *Then the conclusions of Theorem* 3.1 *hold.*

PROOF. Choose the same values for $p_1, \ldots, p_N$ and $q$ as in the proof of Theorem 3.1. Note that, because $\hat{\mathbf{p}} > q^N$ and $\hat{\mathbf{n}}^{\kappa+1} q^{-\mu-N} = o(1)$,

$$(\hat{\mathbf{n}}^{\kappa+1}/\hat{\mathbf{p}})\varphi(q) \leq C\hat{\mathbf{n}}^{\kappa+1} q^{-N} q^{-\mu} = \hat{\mathbf{n}}^{\kappa+1} q^{-\mu-N} \to 0$$

as $\mathbf{n} \Rightarrow \infty$. The end of the proof is entirely similar to that of Theorem 3.1, with Lemma 3.2 instead of Lemma 3.1. □

Analogues of Theorems 3.2–3.4 can also be obtained under assumption (A4″); details are omitted for the sake of brevity.

**4. Numerical results.** In this section, we report the results of a brief Monte Carlo study of the method described in this paper. We mainly consider two models, both in a two-dimensional space ($N = 2$) [writing $(i, j)$ instead of $(i_1, i_2)$ for the sites $\mathbf{i} \in \mathbb{Z}^2$]. For the sake of simplicity, $\mathbf{X}$ (written as $X$) is univariate ($d = 1$).

(a) *Model* 1. Denoting by $\{u_{i,j}, (i,j) \in \mathbb{Z}^2\}$ and $\{e_{i,j}, (i,j) \in \mathbb{Z}^2\}$ two mutually independent i.i.d. $\mathcal{N}(0,1)$ white-noise processes, let

$$Y_{i,j} = g(X_{i,j}) + u_{i,j} \quad \text{with } g(x) := \tfrac{1}{3}e^x + \tfrac{2}{3}e^{-x},$$

where $\{X_{i,j}, (i,j) \in \mathbb{Z}^2\}$ is generated by the spatial autoregression

$$X_{i,j} = \sin(X_{i-1,j} + X_{i,j-1} + X_{i+1,j} + X_{i,j+1}) + e_{i,j}.$$

(b) *Model* 2. Denoting again by $\{e_{i,j}, (i,j) \in \mathbb{Z}^2\}$ an i.i.d. $\mathcal{N}(0,1)$ white-noise process, let $\{Y_{i,j}, (i,j) \in \mathbb{Z}^2\}$ be generated by

$$Y_{i,j} = \sin(Y_{i-1,j} + Y_{i,j-1} + Y_{i+1,j} + Y_{i,j+1}) + e_{i,j},$$

and set

(4.1) $$X_{i,j}^0 := Y_{i-1,j} + Y_{i,j-1} + Y_{i+1,j} + Y_{i,j+1}.$$

Then the prediction function $x \mapsto g(x) := \mathrm{E}[Y_{i,j}|X_{i,j}^0 = x]$ provides the optimal prediction of $Y_{i,j}$ based on $X_{i,j}^0$ in the sense of minimal mean squared prediction error. Note that, in the spatial context, this optimal prediction function $g(\cdot)$ generally differs from the spatial autoregression function itself [here, $\sin(\cdot)$]; see Whittle (1954) for details. Beyond a simple estimation of $g$, we also will investigate the impact, on prediction performance, of including additional spatial lags of $Y_{i,j}$ into the definition of $X_{i,j}$.

Data were simulated from these two models over a rectangular domain of $m \times n$ sites—more precisely, over a grid of the form $\{(i,j)|76 \leq i \leq 75+m, 76 \leq j \leq 75+n\}$, for various values of $m$ and $n$. Each replication was obtained iteratively along the following steps. First, we simulated i.i.d. random variables $e_{ij}$ over the grid $\{(i,j), i = 1, \ldots, 150+m, j = 1, \ldots, 150+n\}$. Next, all initial values of $Y_{ij}$ and $X_{ij}$ being set to zero, we generated $Y_{ij}$'s (or $X_{ij}$'s) over $\{(i,j), i = 1, \ldots, 150+m, j = 1, \ldots, 150+n\}$ recursively, using the spatial autoregressive models. Starting from these generated values, the process was iterated 20 times. The results at the final iteration step for $(i,j)$ inside $\{(i,j)|76 \leq i \leq 75+m, 76 \leq j \leq 75+n\}$ were taken as our simulated $m \times n$ sample. This discarding of peripheral sites allows for a *warming-up zone*, and the first 19 iterations were taken as warming-up steps aiming at achieving stationarity. From the resulting $m \times n$ central data set, we estimated the spatial regression/prediction function using the local linear approach described in this paper. A data-driven choice of the bandwidth in this context would be highly desirable. In view of the lack of theoretical results on this point, we uniformly chose a bandwidth of 0.5 in all our simulations. The simulation results, each with 10 replications, are displayed in Figures 1 and 2 for Models 1 and 2, respectively. Model 1 is a spatial regression model, with the covariates $X_{i,j}$ forming a nonlinear autoregressive process. Inspection of Figure 1 shows that the estimation of the regression function $g(\cdot)$ is quite good and stable, even for sample sizes as small as $m = 10$ and $n = 20$.

Model 2 is a spatial autoregressive model, where $Y_{i,j}$ forms a process with nonlinear spatial autoregression function $\sin(\cdot)$. Various definitions of $X_{i,j}$, involving different spatial lags of $Y_{i,j}$, yield various prediction functions, which are shown in Figures 2(a)–(f). The results in Figures 2(a) and (b) correspond to $X_{i,j} = X_{i,j}^0 := Y_{i-1,j} + Y_{i,j-1} + Y_{i+1,j} + Y_{i,j+1}$, that is, the lags of order $\pm 1$ of $Y_{i,j}$ which also appear in the generating process (4.1). In

$$X_{i,j}^c := Y_{i-2,j} + Y_{i,j-2} + Y_{i-1,j} + Y_{i,j-1} + Y_{i+1,j} + Y_{i,j+1} + Y_{i+2,j} + Y_{i,j+2},$$

thus including lagged values of $Y_{i,j}$ up to order $\pm 2$, in an isotropic way. Nonisotropic choices of $X_{i,j}$ were made in the simulations reported in Figures 2(d)–(f): $X_{i,j}^d := Y_{i-1,j} + Y_{i,j-1}$ in Figure 2(d), $X_{i,j}^e := Y_{i+1,j} + Y_{i,j+1}$ in Figure 2(e) and $X_{i,j}^f := Y_{i-2,j} + Y_{i,j-2} + Y_{i-1,j} + Y_{i,j-1}$ in Figure 2(f).

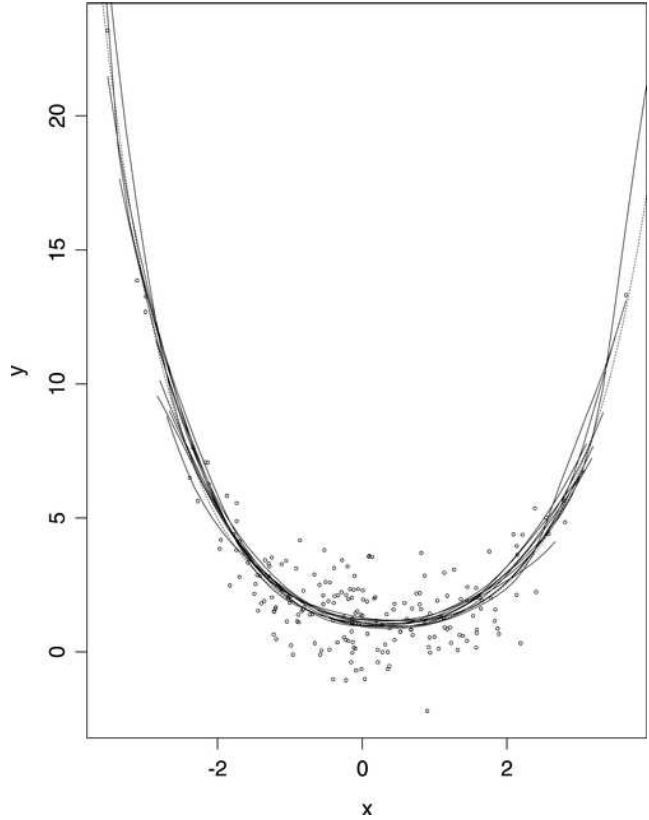

FIG. 1. Simulation for Model 1. *The local linear estimates corresponding to the* 10 *replications* (*solid lines*) *and actual spatial regression curve* (*dotted line*) $g(x) = \mathrm{E}(Y_{ij}|X_{ij} = x) = \frac{1}{3}e^x + \frac{2}{3}e^{-x}$, *for sample size* $m = 10, n = 20$, *with autoregressive spatial covariate* $X_{ij}$. *The scatterplot shows the observations* $(X_{ij}, Y_{ij})$ *corresponding to one typical realization among* 10.

Figure 2(a), the sample sizes $m = 10$ and $n = 20$ are the same as in Figure 1, but the results (still, for 10 replications) are more dispersed. In Figure 2(b), the sample sizes ($m = 30$ and $n = 40$) are slightly larger, and the results (over 10 replications) seem much more stable. These sample sizes therefore were maintained throughout all subsequent simulations. In Figure 2(c), we chose

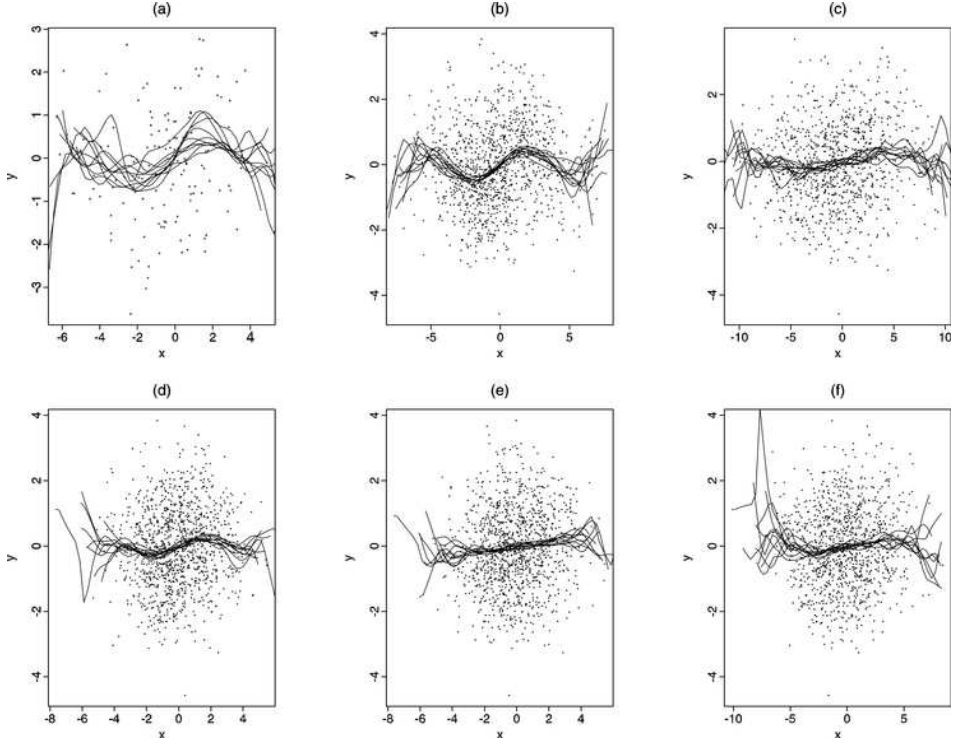

FIG. 2. *Simulation for Model* 2. *The local linear estimates corresponding to the* 10 *replications (solid lines) of the spatial prediction function* $g(x) = \mathrm{E}(Y_{ij}|X_{ij} = x)$, *with sample sizes* $m = 10, n = 20$ *in* (a) *and* $m = 30, n = 40$ *in* (b)–(f), *for different spatial covariates* $X_{ij}$'s: (a) $X_{i,j}^0 := Y_{i-1,j} + Y_{i,j-1} + Y_{i+1,j} + Y_{i,j+1}$; (b) $X_{i,j}^0 := Y_{i-1,j} + Y_{i,j-1} + Y_{i+1,j} + Y_{i,j+1}$; (c) $X_{i,j}^c := Y_{i-2,j} + Y_{i,j-2} + Y_{i-1,j} + Y_{i,j-1} + Y_{i+1,j} + Y_{i,j+1} + Y_{i+2,j} + Y_{i,j+2}$; (d) $X_{i,j}^d := Y_{i-1,j} + Y_{i,j-1}$; (e) $X_{i,j}^e := Y_{i+1,j} + Y_{i,j+1}$; *and* (f) $X_{i,j}^f := Y_{i-2,j} + Y_{i,j-2} + Y_{i-1,j} + Y_{i,j-1}$. *The scatterplot shows the observations* $(X_{ij}, Y_{ij})$ *corresponding to one typical realization among* 10.

A more systematic simulation study certainly would be welcome. However, it seems that, even in very small samples (see Figure 1), the performance of our method is excellent in pure spatial regression problems (with spatially correlated covariates), while larger samples are required in spatial *auto*regression models. This difference is probably strongly related to differences in the corresponding *noise-to-signal ratios*. Letting $g(x) = \mathrm{E}(Y|X=x)$ and $\varepsilon = Y - g(X)$, the noise-to-signal ratio is defined as $\mathrm{Var}(\varepsilon)/\mathrm{Var}(g(X))$; see, for example, Chapter 4 in Fan and Gijbels (1996) for details. In a classical regression setting, independence is generally assumed between $X$ and $\varepsilon$, so that this ratio, in simulations, can be set in advance. Such an independence assumption cannot be made in a spatial series context, but empirical versions of the ratio nevertheless can be computed from each replication, then averaged, providing estimated values. In Model 1 this estimated value (averaged over the 10 replications) of the noise-to-signal ratio is 0.214. The values for the six versions of Model 2 (still, averaged over 10 replications) are much larger: (a) 12.037, (b) 13.596, (c) 43.946, (d) 47.442, (e) 116.334 and (f) 88.287.

## 5. Proofs.

5.1. *Proof of Lemma* 2.2. The proof of Lemma 2.2 relies on two intermediate results. The first one is a lemma borrowed from Ibragimov and Linnik (1971) or Deo (1973), to which we refer for a proof.

LEMMA 5.1. (i) *Suppose that* (A1) *holds. Let* $\mathcal{L}_r(\mathcal{F})$ *denote the class of $\mathcal{F}$-measurable random variables $\xi$ satisfying* $\|\xi\|_r := (\mathrm{E}|\xi|^r)^{1/r} < \infty$. *Let* $\mathbf{X} \in \mathcal{L}_r(\mathcal{B}(\mathcal{S}))$ *and* $Y \in \mathcal{L}_s(\mathcal{B}(\mathcal{S}'))$. *Then for any* $1 \leq r, s, h < \infty$ *such that* $r^{-1} + s^{-1} + h^{-1} = 1$,

$$(5.1) \qquad |\mathrm{E}[\mathbf{X}Y] - \mathrm{E}[\mathbf{X}]\mathrm{E}[Y]| \leq C\|\mathbf{X}\|_r\|Y\|_s[\alpha(\mathcal{S}, \mathcal{S}')]^{1/h},$$

*where* $\|X\|_2 := \|(X'X)^{1/2}\|_r$.

(ii) *If, moreover,* $\|\mathbf{X}\| := (\mathbf{X}^r\mathbf{X})^{1/2}$ *and* $|Y|$ *are* P-*a.s. bounded, the right-hand side of* (5.1) *can be replaced by* $C\alpha(\mathcal{S}, \mathcal{S}')$.

The second one is a lemma of independent interest, which plays a crucial role here and in the subsequent sections. For the sake of generality, and in order for this lemma to apply beyond the specific context of this paper, we do not necessarily assume that the mixing coefficient $\alpha$ takes the form imposed in assumption (A4).

Before stating the lemma, let us first introduce some further notation. Let

$$A_{\mathbf{n}} = (\hat{\mathbf{n}}b_{\mathbf{n}}^d)^{-1/2} \sum_{\mathbf{j} \in \mathcal{I}_{\mathbf{n}}} \eta_{\mathbf{j}}(x)$$

and

$$\mathrm{Var}(A_{\mathbf{n}}) = (\hat{\mathbf{n}}b_{\mathbf{n}}^d)^{-1} \sum_{\mathbf{j} \in \mathcal{I}_{\mathbf{n}}} \mathrm{E}[\Delta_{\mathbf{j}}^2(\mathbf{x})] + (\hat{\mathbf{n}}b_{\mathbf{n}}^d)^{-1} \sum_{\{\mathbf{i},\mathbf{j} \in \mathcal{I}_{\mathbf{n}} | \exists k : i_k \neq j_k\}} \mathrm{E}[\Delta_{\mathbf{i}}(\mathbf{x})\Delta_{\mathbf{j}}(\mathbf{x})]$$
$$:= \tilde{I}(\mathbf{x}) + \tilde{R}(\mathbf{x}), \quad \text{say},$$

where $\eta_{\mathbf{j}}(\mathbf{x}) := Z_{\mathbf{j}} K_{\mathbf{c}}(\mathbf{x} - \mathbf{X}_{\mathbf{j}})$ and $\Delta_{\mathbf{j}}(\mathbf{x}) := \eta_{\mathbf{j}}(\mathbf{x}) - \mathrm{E}\eta_{\mathbf{j}}(\mathbf{x})$. For any $\mathbf{c_n} := (c_{\mathbf{n}1}, \ldots, c_{\mathbf{n}N}) \in \mathbb{Z}^N$ with $1 < c_{\mathbf{n}k} < n_k$ for all $k = 1, \ldots, N$, define $\tilde{J}_1(\mathbf{x}) :=$

$b_{\mathbf{n}}^{\delta d/(4+\delta)+d} \prod_{k=1}^{N}(n_k c_{\mathbf{n}k})$ and

$$\tilde{J}_2(\mathbf{x}) := b_{\mathbf{n}}^{2d/(2+\delta)} \hat{\mathbf{n}} \sum_{k=1}^{N} \left( \sum_{\substack{|j_s|=1 \\ s=1,\ldots,k-1}}^{n_s} \sum_{|j_k|=c_{\mathbf{n}k}}^{n_k} \sum_{\substack{|j_s|=1 \\ s=k+1,\ldots,N}}^{n_s} \{\varphi(j_1,\ldots,j_N)\}^{\delta/(2+\delta)} \right).$$

LEMMA 5.2. *Let $\{(Y_{\mathbf{j}}, \mathbf{X}_{\mathbf{j}}); \mathbf{j} \in \mathbb{Z}^N\}$ denote a stationary spatial process with general mixing coefficient*

$$\begin{aligned}\varphi(\mathbf{j}) &= \varphi(j_1,\ldots,j_N) \\ &:= \sup\{|\mathrm{P}(AB) - \mathrm{P}(A)\mathrm{P}(B)| : A \in \mathcal{B}(\{Y_{\mathbf{i}}, \mathbf{X}_{\mathbf{i}}\}), B \in \mathcal{B}(\{Y_{\mathbf{i}+\mathbf{j}}, \mathbf{X}_{\mathbf{i}+\mathbf{j}}\})\},\end{aligned}$$

*and assume that assumptions* (A1), (A2) *and* (A5) *hold. Then*

(5.2) $$|\tilde{R}(\mathbf{x})| \leq C(\hat{\mathbf{n}} b_{\mathbf{n}}^d)^{-1}[\tilde{J}_1(\mathbf{x}) + \tilde{J}_2(\mathbf{x})].$$

*If furthermore $\varphi(j_1,\ldots,j_N)$ takes the form $\varphi(\|\mathbf{j}\|)$, then*

(5.3) $$\tilde{J}_2(\mathbf{x}) \leq C b_{\mathbf{n}}^{2d/(2+\delta)} \hat{\mathbf{n}} \sum_{k=1}^{N} \left( \sum_{t=c_{\mathbf{n}k}}^{\|\mathbf{n}\|} t^{N-1} \{\varphi(t)\}^{\delta/(2+\delta)} \right).$$

PROOF. Set $L = L_{\mathbf{n}} = b_{\mathbf{n}}^{-2d/(4+\delta)}$. Defining $Z_{1\mathbf{j}} := Z_{\mathbf{j}} I_{\{|Z_{\mathbf{j}}| \leq L\}}$ and $Z_{2\mathbf{j}} := Z_{\mathbf{j}} I_{\{|Z_{\mathbf{j}}| > L\}}$, let

$$\eta_{i\mathbf{j}}(\mathbf{x}) := Z_{i\mathbf{j}} K_{\mathbf{c}}(\mathbf{x} - \mathbf{X}_{\mathbf{j}}) \quad \text{and} \quad \Delta_{i\mathbf{j}}(\mathbf{x}) := \eta_{i\mathbf{j}}(\mathbf{x}) - E\eta_{i\mathbf{j}}(\mathbf{x}), \qquad i=1,2.$$

Then $Z_{\mathbf{j}} = Z_{1\mathbf{j}} + Z_{2\mathbf{j}}, \Delta_{\mathbf{j}}(\mathbf{x}) = \Delta_{1\mathbf{j}}(\mathbf{x}) + \Delta_{2\mathbf{j}}(\mathbf{x})$, and hence

(5.4) $$\begin{aligned}E\Delta_{\mathbf{j}}(\mathbf{x})\Delta_{\mathbf{i}}(\mathbf{x}) = &E\Delta_{1\mathbf{j}}(\mathbf{x})\Delta_{1\mathbf{i}}(\mathbf{x}) + E\Delta_{1\mathbf{j}}(\mathbf{x})\Delta_{2\mathbf{i}}(\mathbf{x}) \\ &+ E\Delta_{2\mathbf{j}}(\mathbf{x})\Delta_{1\mathbf{i}}(\mathbf{x}) + E\Delta_{2\mathbf{j}}(\mathbf{x})\Delta_{2\mathbf{i}}(\mathbf{x}).\end{aligned}$$

First, we note that

$$\begin{aligned}b_{\mathbf{n}}^{-d} &|E\Delta_{1\mathbf{j}}(\mathbf{x})\Delta_{2\mathbf{i}}(\mathbf{x})| \\ &\leq \{b_{\mathbf{n}}^{-d} E\eta_{1\mathbf{j}}^2(\mathbf{x})\}^{1/2} \{b_{\mathbf{n}}^{-d} E\eta_{2\mathbf{i}}^2(\mathbf{x})\}^{1/2} \\ &\leq \{b_{\mathbf{n}}^{-d} E Z_{1\mathbf{j}}^2 K_{\mathbf{c}}^2((\mathbf{x}-\mathbf{X}_{\mathbf{j}})/b_{\mathbf{n}})\}^{1/2} \{b_{\mathbf{n}}^{-d} E Z_{2\mathbf{i}}^2 K_{\mathbf{c}}^2((\mathbf{x}-\mathbf{X}_{\mathbf{j}})/b_{\mathbf{n}})\}^{1/2} \\ &\leq C\{b_{\mathbf{n}}^{-d} E|Z_{\mathbf{i}}|^2 I_{\{|Z_{\mathbf{i}}|>L\}} K_{\mathbf{c}}((\mathbf{x}-\mathbf{X}_{\mathbf{1}})/b_{\mathbf{n}})\}^{1/2} \\ &\leq C\{L^{-\delta} b_{\mathbf{n}}^{-d} E|Z_{\mathbf{j}}|^{2+\delta} I_{\{|Z_{\mathbf{j}}|>L\}} K_{\mathbf{c}}((\mathbf{x}-\mathbf{X}_{\mathbf{1}})/b_{\mathbf{n}})\}^{1/2} \\ &\leq C L_{\mathbf{n}}^{-\delta/2} = C b_{\mathbf{n}}^{\delta d/(4+\delta)}.\end{aligned}$$

Similarly,
$$b_{\mathbf{n}}^{-d}|\mathrm{E}\Delta_{2\mathbf{j}}(\mathbf{x})\Delta_{1\mathbf{i}}(\mathbf{x})| \leq CL_{\mathbf{n}}^{-\delta/2} = Cb_{\mathbf{n}}^{\delta d/(4+\delta)} \quad \text{and}$$
$$b_{\mathbf{n}}^{-d}|\mathrm{E}\Delta_{2\mathbf{j}}(\mathbf{x})\Delta_{2\mathbf{i}}(\mathbf{x})| \leq Cb_{\mathbf{n}}^{2\delta d/(4+\delta)}.$$

Next, for $\mathbf{i} \neq \mathbf{j}$, letting $K_{\mathbf{n}}(\mathbf{x}) := (1/b_{\mathbf{n}}^d)K(\mathbf{x}/b_{\mathbf{n}})$ and $K_{\mathbf{cn}}(\mathbf{x}) := (1/b_{\mathbf{n}}^d)K_{\mathbf{c}}(\mathbf{x}/b_{\mathbf{n}})$,

$$b_{\mathbf{n}}^{-d}\mathrm{E}\Delta_{1\mathbf{j}}(\mathbf{x})\Delta_{1\mathbf{i}}(\mathbf{x})$$
$$= b_{\mathbf{n}}^d\{\mathrm{E}Z_{1\mathbf{i}}Z_{1\mathbf{j}}K_{\mathbf{cn}}(\mathbf{x}-\mathbf{X_i})K_{\mathbf{cn}}(\mathbf{x}-\mathbf{X_j})$$
$$\qquad - \mathrm{E}Z_{1\mathbf{i}}K_{\mathbf{cn}}(\mathbf{x}-\mathbf{X_i})\mathrm{E}Z_{1\mathbf{j}}K_{\mathbf{cn}}(\mathbf{x}-\mathbf{X_j})\}$$
$$= b_{\mathbf{n}}^d \iint K_{\mathbf{cn}}(\mathbf{x}-\mathbf{u})K_{\mathbf{cn}}(\mathbf{x}-\mathbf{v})$$
$$\qquad \times \{g_{1\mathbf{ij}}(\mathbf{u},\mathbf{v})f_{\mathbf{i},\mathbf{j}}(\mathbf{u},\mathbf{v}) - g_1^{(1)}(\mathbf{u})g_1^{(1)}(\mathbf{v})f(\mathbf{u})f(\mathbf{v})\}\,d\mathbf{u}\,d\mathbf{v},$$

where $g_{1\mathbf{ij}}(\mathbf{u},\mathbf{v}) := \mathrm{E}(Z_{1\mathbf{i}}Z_{1\mathbf{j}}|\mathbf{X_i}=\mathbf{u},\mathbf{X_j}=\mathbf{v})$, and $g_1^{(1)}(\mathbf{u}) := \mathrm{E}(Z_{1\mathbf{i}}|\mathbf{X_i}=\mathbf{u})$. Since, by definition, $|Z_{1\mathbf{i}}| \leq L_{\mathbf{n}}$, we have that $|g_{1\mathbf{ij}}(\mathbf{u},\mathbf{v})| \leq L_{\mathbf{n}}^2$ and $|g_1^{(1)}(\mathbf{u}) \times g_1^{(1)}(\mathbf{v})| \leq L_{\mathbf{n}}^2$. Thus

$$|g_{1\mathbf{ij}}(\mathbf{u},\mathbf{v})f_{\mathbf{i},\mathbf{j}}(\mathbf{u},\mathbf{v}) - g_1^{(1)}(\mathbf{u})g_1^{(1)}(\mathbf{v})f(\mathbf{u})f(\mathbf{v})|$$
$$\leq |g_{1\mathbf{ij}}(\mathbf{u},\mathbf{v})(f_{\mathbf{i},\mathbf{j}}(\mathbf{u},\mathbf{v}) - f(\mathbf{u})f(\mathbf{v}))|$$
$$\quad + |(g_{1\mathbf{ij}}(\mathbf{u},\mathbf{v}) - g_1^{(1)}(\mathbf{u})g_1^{(1)}(\mathbf{v}))f(\mathbf{u})f(\mathbf{v})|$$
$$\leq L_{\mathbf{n}}^2|f_{\mathbf{i},\mathbf{j}}(\mathbf{u},\mathbf{v}) - f(\mathbf{u})f(\mathbf{v})| + 2L_{\mathbf{n}}^2 f(\mathbf{u})f(\mathbf{v}).$$

It then follows from (A1) and the Lebesgue density theorem [see Chapter 2 of Devroye and Györfi (1985)] that

$$b_{\mathbf{n}}^{-d}|\mathrm{E}\Delta_{1\mathbf{j}}(\mathbf{x})\Delta_{1\mathbf{i}}(\mathbf{x})|$$
$$\leq b_{\mathbf{n}}^d \iint K_{\mathbf{cn}}(\mathbf{x}-\mathbf{u})K_{\mathbf{cn}}(\mathbf{x}-\mathbf{v})L_{\mathbf{n}}^2|f_{\mathbf{i},\mathbf{j}}(\mathbf{u},\mathbf{v}) - f(\mathbf{u})f(\mathbf{v})|\,d\mathbf{u}\,d\mathbf{v}$$
$$\quad + b_{\mathbf{n}}^d \iint 2L_{\mathbf{n}}^2 f(\mathbf{u})f(\mathbf{v})\,d\mathbf{u}\,d\mathbf{v}$$
$$\leq Cb_{\mathbf{n}}^d\left(L_{\mathbf{n}}^2\left\{\int K_{\mathbf{cn}}(\mathbf{x}-\mathbf{u})\,d\mathbf{u}\right\}^2 + 2L_{\mathbf{n}}^2\left\{\int K_{\mathbf{n}}(\mathbf{x}-\mathbf{u})f(\mathbf{u})\,d\mathbf{u}\right\}^2\right)$$
$$\leq Cb_{\mathbf{n}}^d L_{\mathbf{n}}^2 = Cb_{\mathbf{n}}^{\delta d/(4+\delta)}.$$

(5.5)

Thus, by (5.4) and (5.5),

(5.6) $\qquad b_{\mathbf{n}}^{-d}|\mathrm{E}\Delta_{\mathbf{j}}(\mathbf{x})\Delta_{\mathbf{i}}(\mathbf{x})| \leq CL_{\mathbf{n}}^{-\delta/2} + Cb_{\mathbf{n}}^d L_{\mathbf{n}}^2 = Cb_{\mathbf{n}}^{\delta d/(4+\delta)}.$

Let $\mathbf{c_n} = (c_{\mathbf{n}1},\ldots,c_{\mathbf{n}N}) \in \mathbb{R}^N$ be a sequence of vectors with positive components. Define

$$\mathcal{S}_1 := \{\mathbf{i} \neq \mathbf{j} \in \mathcal{I}_{\mathbf{n}} : |j_k - i_k| \leq c_{\mathbf{n}k} \text{ for all } k = 1,\ldots,N\}$$

and

$$S_2 := \{\mathbf{i}, \mathbf{j} \in \mathcal{I}_\mathbf{n} : |j_k - i_k| > c_{\mathbf{n}k} \text{ for some } k = 1, \ldots, N\}.$$

Clearly, $\text{Card}(S_1) \leq 2^N \hat{\mathbf{n}} \prod_{k=1}^N c_{\mathbf{n}k}$. Splitting $\tilde{R}(\mathbf{x})$ into $(\hat{\mathbf{n}} b_\mathbf{n}^d)^{-1}(J_1 + J_2)$, with

$$J_\ell := \sum\sum_{\mathbf{i},\mathbf{j}\in S_\ell} \mathrm{E}\Delta_\mathbf{j}(\mathbf{x})\Delta_\mathbf{i}(\mathbf{x}), \qquad \ell = 1, 2,$$

it follows from (5.6) that

(5.7) $$|J_1| \leq C b_\mathbf{n}^{\delta d/(4+\delta)+d} \text{Card}(S_1) \leq 2^N C b_\mathbf{n}^{\delta d/(4+\delta)+d} \hat{\mathbf{n}} \prod_{k=1}^N c_{\mathbf{n}k}.$$

Turning to $J_2$, we have $|J_2| \leq \sum\sum_{\mathbf{i},\mathbf{j}\in S_2} |\mathrm{E}\Delta_\mathbf{j}(\mathbf{x})\Delta_\mathbf{i}(\mathbf{x})|$. Lemma 5.1, with $r = s = 2 + \delta$ and $h = (2+\delta)/\delta$, yields

$$\begin{aligned}|\mathrm{E}\Delta_\mathbf{j}(\mathbf{x})\Delta_\mathbf{i}(\mathbf{x})| &\leq C(\mathrm{E}|Z_\mathbf{i} K_\mathbf{c}((\mathbf{x}-\mathbf{X_i})/b_\mathbf{n})|^{2+\delta})^{2/(2+\delta)}\{\varphi(\mathbf{j}-\mathbf{i})\}^{\delta/(2+\delta)} \\ &\leq C b_\mathbf{n}^{2d/(2+\delta)} (b_\mathbf{n}^{-d}\mathrm{E}|Z_\mathbf{i} K_\mathbf{c}((\mathbf{x}-\mathbf{X_i})/b_\mathbf{n})|^{2+\delta})^{2/(2+\delta)}\{\varphi(\mathbf{j}-\mathbf{i})\}^{\delta/(2+\delta)} \\ &\leq C b_\mathbf{n}^{2d/(2+\delta)}\{\varphi(\mathbf{j}-\mathbf{i})\}^{\delta/(2+\delta)}.\end{aligned}$$
(5.8)

Hence,

(5.9) $$|J_2| \leq C b_\mathbf{n}^{2d/(2+\delta)} \sum\sum_{\mathbf{i},\mathbf{j}\in S_2}\{\varphi(\mathbf{j}-\mathbf{i})\}^{\delta/(2+\delta)} := C b_\mathbf{n}^{2d/(2+\delta)} \Sigma_2, \qquad \text{say}.$$

We now analyze the quantity $\Sigma_2$ in detail. For any $N$-tuple $\mathbf{0} \neq \boldsymbol{\ell} = (\ell_1,\ldots,\ell_N) \in \{0,1\}^N$, set

$$\begin{aligned}\mathcal{S}(\ell_1,\ldots,\ell_N) := \{\mathbf{i},\mathbf{j} \in I_\mathbf{n} : &|j_k - i_k| > c_{\mathbf{n}k} \text{ if } \ell_k = 1 \text{ and} \\ &|j_k - i_k| \leq c_{\mathbf{n}k} \text{ if } \ell_k = 0, k = 1,\ldots,N\}\end{aligned}$$

and

$$V(\ell_1,\ldots,\ell_N) := \sum\sum_{\mathbf{i},\mathbf{j}\in \mathcal{S}(\ell_1,\ldots,\ell_N)} \{\varphi(\mathbf{j}-\mathbf{i})\}^{\delta/(2+\delta)}.$$

Then

(5.10) $$\Sigma_2 = \sum\sum_{\mathbf{i},\mathbf{j}\in S_2} \{\varphi(\mathbf{j}-\mathbf{i})\}^{\delta/(2+\delta)} = \sum_{\mathbf{0}\neq\boldsymbol{\ell}\in\{0,1\}^N} V(\ell_1,\ldots,\ell_N).$$

Without loss of generality, consider $V(1,0,\ldots,0)$. Because $\sum_{|i_k-j_k|>c_{\mathbf{n}k}}(\cdots)$ decomposes into $\sum_{i_k=1}^{n_k-c_{\mathbf{n}k}-1}\sum_{j_k=i_k+c_{\mathbf{n}k}+1}^{n_k}(\cdots) + \sum_{j_k=1}^{n_k-c_{\mathbf{n}k}-1}\sum_{i_k=j_k+c_{\mathbf{n}k}+1}^{n_k}(\cdots),$

and $\sum_{|i_k-j_k|\leq c_{\mathbf{n}k}}(\cdots)$ into $\sum_{i_k=1}^{n_k-c_{\mathbf{n}k}}\sum_{j_k=i_k+1}^{i_k+c_{\mathbf{n}k}}(\cdots) + \sum_{j_k=1}^{n_k-c_{\mathbf{n}k}}\sum_{i_k=j_k+1}^{j_k+c_{\mathbf{n}k}}(\cdots)$, we have

$$V(1,0,\ldots,0)$$
$$= \sum_{|i_1-j_1|>c_{\mathbf{n}1}} \sum_{|i_2-j_2|\leq c_{\mathbf{n}2}} \cdots \sum_{|i_N-j_N|\leq c_{\mathbf{n}N}} \{\varphi(j_1-i_1,\ldots,j_N-i_N)\}^{\delta/(2+\delta)}$$
$$\leq \hat{\mathbf{n}} \left\{\sum_{j_1=c_{\mathbf{n}1}}^{n_1} + \sum_{-j_1=c_{\mathbf{n}1}}^{n_1}\right\}\left\{\sum_{j_2=1}^{c_{\mathbf{n}2}} + \sum_{-j_2=1}^{c_{\mathbf{n}2}}\right\}\cdots$$
$$\left\{\sum_{j_N=1}^{c_{\mathbf{n}N}} + \sum_{-j_N=1}^{c_{\mathbf{n}N}}\right\}\{\varphi(j_1,\ldots,j_N)\}^{\delta/(2+\delta)}$$
$$\leq \hat{\mathbf{n}} \sum_{|j_1|=c_{\mathbf{n}1}}^{n_1} \sum_{|j_2|=1}^{c_{\mathbf{n}2}} \cdots \sum_{|j_N|=1}^{c_{\mathbf{n}N}} \{\varphi(j_1,\ldots,j_N)\}^{\delta/(2+\delta)}$$
$$\leq \hat{\mathbf{n}} \sum_{|j_1|=c_{\mathbf{n}1}}^{n_1} \sum_{|j_2|=1}^{n_2} \cdots \sum_{|j_N|=1}^{n_N} \{\varphi(j_1,\ldots,j_N)\}^{\delta/(2+\delta)}.$$

More generally,

$$(5.11) \quad V(\ell_1,\ell_2,\ldots,\ell_N) \leq \hat{\mathbf{n}} \sum_{|j_1|}\cdots\sum_{|j_k|}\cdots\sum_{|j_N|}\{\varphi(j_1,\ldots,j_N)\}^{\delta/(2+\delta)},$$

where the sums $\sum_{|j_k|}$ run over all values of $j_k$ such that $1 \leq |j_k| \leq n_k$ if $\ell_k = 0$, and such that $c_{\mathbf{n}1} \leq |j_k| \leq n_k$ if $\ell_k = 1$. Since the summands are nonnegative, for $1 \leq c_{\mathbf{n}k} \leq n_k$, we have $\sum_{|j_k|=c_{\mathbf{n}k}}^{n_k}(\cdots) \leq \sum_{|j_k|=1}^{n_k}(\cdots)$, and (5.9)–(5.11) imply

$$|J_2| \leq Cb_{\mathbf{n}}^{2d/(2+\delta)}\hat{\mathbf{n}}$$
$$(5.12) \quad \times \sum_{k=1}^{N}\left(\sum_{|j_1|=1}^{n_1}\cdots\sum_{|j_{k-1}|=1}^{n_{k-1}}\sum_{|j_k|=c_{\mathbf{n}k}}^{n_k}\sum_{|j_{k+1}|=1}^{n_{k+1}}\cdots\right.$$
$$\left.\sum_{|j_N|=1}^{n_N}\{\varphi(j_1,\ldots,j_N)\}^{\delta/(2+\delta)}\right).$$

Thus, (5.2) is a consequence of (5.7) and (5.12). If, furthermore, $\varphi(j_1,\ldots,j_N)$ depends on $\|\mathbf{j}\|$ only, then

$$\sum_{|j_1|=1}^{n_1}\cdots\sum_{|j_{k-1}|=1}^{n_{k-1}}\sum_{|j_k|=c_{\mathbf{n}k}}^{n_k}\sum_{|j_{k+1}|=1}^{n_{k+1}}\cdots\sum_{|j_N|=1}^{n_N}\{\varphi(\|\mathbf{j}\|)\}^{\delta/(2+\delta)}$$
$$\leq \sum_{|j_1|=1}^{n_1}\cdots\sum_{|j_{k-1}|=1}^{n_{k-1}}\sum_{|j_k|=c_{\mathbf{n}k}}^{n_k}\sum_{|j_{k+1}|=1}^{n_{k+1}}\cdots\sum_{|j_{N-1}|=1}^{n_{N-1}}\sum_{t^2=j_1^2+\cdots+j_{N-1}^2+1}^{j_1^2+\cdots+j_{N-1}^2+n_N^2}\{\varphi(t)\}^{\delta/(2+\delta)}$$

$$\leq \sum_{t=c_{\mathbf{n}k}}^{\|\mathbf{n}\|} \sum_{|j_1|=1}^{t} \cdots \sum_{|j_{N-1}|=1}^{t} \{\varphi(t)\}^{\delta/(2+\delta)} \leq \sum_{t=c_{\mathbf{n}k}}^{\|\mathbf{n}\|} t^{N-1} \{\varphi(t)\}^{\delta/(2+\delta)};$$

(5.3) follows. □

PROOF OF LEMMA 2.2. Observe that

(5.13)
$$\tilde{I}(\mathbf{x}) = b_{\mathbf{n}}^{-d} \mathrm{E} \Delta_{\mathbf{j}}^2(\mathbf{x}) = b_{\mathbf{n}}^{-d} [\mathrm{E} \eta_{\mathbf{j}}^2 - (\mathrm{E} \eta_{\mathbf{j}})^2]$$
$$= b_{\mathbf{n}}^{-d} [\mathrm{E} Z_{\mathbf{j}}^2 K_{\mathbf{c}}^2((\mathbf{x}-\mathbf{X_j})/b_{\mathbf{n}}) - \{\mathrm{E} Z_{\mathbf{j}} K_{\mathbf{c}}((\mathbf{x}-\mathbf{X_j})/b_{\mathbf{n}})\}^2].$$

Under assumption (A5), by the Lebesgue density theorem,

$$\lim_{\mathbf{n}\to\infty} \int_{\mathbb{R}^d} b_{\mathbf{n}}^{-d} \mathrm{E}[Z_{\mathbf{j}}^2 | \mathbf{X_j} = \mathbf{u}] K_{\mathbf{c}}^2((\mathbf{x}-\mathbf{u})/b_{\mathbf{n}}) f(\mathbf{u}) \, d\mathbf{u} = g^{(2)}(\mathbf{x}) f(\mathbf{x}) \int_{\mathbb{R}^d} K_{\mathbf{c}}^2(\mathbf{u}) \, d\mathbf{u},$$

$$\lim_{\mathbf{n}\to\infty} \int_{\mathbb{R}^d} b_{\mathbf{n}}^{-d} \mathrm{E}[Z_{\mathbf{j}} | \mathbf{X_j} = \mathbf{u}] K_{\mathbf{c}}((\mathbf{x}-\mathbf{u})/b_{\mathbf{n}}) f(\mathbf{u}) \, d\mathbf{u} = g^{(1)}(\mathbf{x}) f(\mathbf{x}) \int_{\mathbb{R}^d} K(\mathbf{u}) \, d\mathbf{u},$$

where $g^{(i)}(\mathbf{x}) := \mathrm{E}[Z_{\mathbf{j}}^i | \mathbf{X_j} = \mathbf{x}]$ for $i = 1, 2$. It is easily seen that $b_{\mathbf{n}}^{-d} \{\mathrm{E} Z_{\mathbf{j}} K_{\mathbf{c}}((\mathbf{x}-\mathbf{X_j})/b_{\mathbf{n}})\}^2 \to 0$. Thus, from (5.13),

(5.14)
$$\lim_{\mathbf{n}\to\infty} \tilde{I}(\mathbf{x}) = g^{(2)}(\mathbf{x}) f(\mathbf{x}) \int_{\mathbb{R}^d} K_{\mathbf{c}}^2(\mathbf{u}) \, d\mathbf{u},$$

where $g^{(2)}(\mathbf{x}) = \mathrm{E}\{Z_{\mathbf{j}}^2 | \mathbf{X_j} = \mathbf{x}\} = \mathrm{E}\{(Y_{\mathbf{j}} - g(\mathbf{x}))^2 | \mathbf{X_j} = \mathbf{x}\} = \mathrm{Var}\{Y_{\mathbf{j}} | \mathbf{X_j} = \mathbf{x}\}$.

Let $c_{\mathbf{n}k}^a := b_{\mathbf{n}}^{-\delta d/(2+\delta)} \to \infty$. Clearly, $c_{\mathbf{n}k} < n_k$ because $n_k b_{\mathbf{n}}^{\delta d/(2+\delta)a} > 1$ for all $k$. Apply Lemma 5.2. Since, due to the fact that $a > (4+\delta)N/(2+\delta)$, and $N/(2+\delta)a < 1/(4+\delta)$

(5.15)
$$(\hat{\mathbf{n}} b_{\mathbf{n}}^d)^{-1} \tilde{J}_2 \leq C \sum_{k=1}^{N} \left( c_{\mathbf{n}k}^a \sum_{t=c_{\mathbf{n}k}}^{\infty} t^{N-1} \{\varphi(t)\}^{\delta/(2+\delta)} \right) \to 0$$

because $c_{\mathbf{n}k} \to \infty$, (5.3) and assumption (A4) imply that

$$(\hat{\mathbf{n}} b_{\mathbf{n}}^d)^{-1} \tilde{J}_1 \leq C b_{\mathbf{n}}^{\delta d/(4+\delta)} c_{\mathbf{n}1} \cdots c_{\mathbf{n}N} = C b_{\mathbf{n}}^{\delta d/(4+\delta)} b_{\mathbf{n}}^{-\delta dN/(2+\delta)a} \to 0,$$

hence, by (5.2), that

(5.16) $\qquad |\tilde{R}(x)| = (\hat{\mathbf{n}} b_{\mathbf{n}}^d)^{-1} |\tilde{J}(x)| \leq C(\hat{\mathbf{n}} b_{\mathbf{n}}^d)^{-1} (\tilde{J}_1 + \tilde{J}_2) \to 0.$

Finally, (2.7) follows from (5.14) and (5.16), which completes the proof of Lemma 2.2. □

PROOF OF LEMMA 2.3. From (2.5) and the definition of $A_{\mathbf{n}}$ [recall that $a_0 = g(\mathbf{x})$, $\mathbf{a}_1 = g'(\mathbf{x})$],

$$\mathrm{E}[A_{\mathbf{n}}] = (\hat{\mathbf{n}} b_{\mathbf{n}}^d)^{1/2} b_{\mathbf{n}}^{-d} \mathrm{E}[Z_{\mathbf{j}}] K_{\mathbf{c}} \left( \frac{\mathbf{X_j} - \mathbf{x}}{b_{\mathbf{n}}} \right)$$

$$= (\hat{\mathbf{n}} b_{\mathbf{n}}^d)^{1/2} \, b_{\mathbf{n}}^{-d} \mathrm{E}(Y_{\mathbf{j}} - a_0 - \mathbf{a}_1^\tau(\mathbf{X}_{\mathbf{j}} - \mathbf{x})) K_{\mathbf{c}}\!\left(\frac{\mathbf{X}_{\mathbf{j}} - \mathbf{x}}{b_{\mathbf{n}}}\right)$$

$$= (\hat{\mathbf{n}} b_{\mathbf{n}}^d)^{1/2} b_{\mathbf{n}}^{-d} \mathrm{E}(g(\mathbf{X}_{\mathbf{j}}) - a_0 - \mathbf{a}_1^\tau(\mathbf{X}_{\mathbf{j}} - \mathbf{x})) K_{\mathbf{c}}\!\left(\frac{\mathbf{X}_{\mathbf{j}} - \mathbf{x}}{b_{\mathbf{n}}}\right)$$

$$= (\hat{\mathbf{n}} b_{\mathbf{n}}^d)^{1/2} b_{\mathbf{n}}^{-d} \mathrm{E}(\mathbf{X}_{\mathbf{j}} - \mathbf{x})^\tau$$
$$\times g''(\mathbf{x} + \boldsymbol{\xi}(\mathbf{X}_{\mathbf{j}} - \mathbf{x}))(\mathbf{X}_{\mathbf{j}} - \mathbf{x}) K_{\mathbf{c}}\!\left(\frac{\mathbf{X}_{\mathbf{j}} - \mathbf{x}}{b_{\mathbf{n}}}\right) \quad \text{(where } |\boldsymbol{\xi}| < 1\text{)}$$

$$= (\hat{\mathbf{n}} b_{\mathbf{n}}^d)^{1/2} b_{\mathbf{n}}^2 \, b_{\mathbf{n}}^{-d} \operatorname{tr} \mathrm{E}\!\left[g''(\mathbf{x} + \boldsymbol{\xi}(\mathbf{X}_{\mathbf{j}} - \mathbf{x})) \frac{\mathbf{X}_{\mathbf{j}} - \mathbf{x}}{b_{\mathbf{n}}}\!\left(\frac{\mathbf{X}_{\mathbf{j}} - \mathbf{x}}{b_{\mathbf{n}}}\right)^\tau\right] K_{\mathbf{c}}\!\left(\frac{\mathbf{X}_{\mathbf{j}} - \mathbf{x}}{b_{\mathbf{n}}}\right);$$

the lemma follows via assumption (A3). □

PROOF OF LEMMA 3.1. The proof consists of two parts and an additional lemma (Lemma 5.3). Recalling that

(5.17) $\quad \eta_{\mathbf{j}}(\mathbf{x}) := Z_{\mathbf{j}} K_{\mathbf{c}}(\mathbf{x} - \mathbf{X}_{\mathbf{j}}) \quad \text{and} \quad \Delta_{\mathbf{j}}(\mathbf{x}) := \eta_{\mathbf{j}}(\mathbf{x}) - \mathrm{E}\eta_{\mathbf{j}}(\mathbf{x}),$

define $\zeta_{\mathbf{nj}} := b_{\mathbf{n}}^{-d/2} \Delta_{\mathbf{j}}$, and let $S_{\mathbf{n}} := \sum_{j_k=1; k=1,\ldots,N}^{n_k} \zeta_{\mathbf{nj}}$. Then

$$\hat{\mathbf{n}}^{-1/2} S_{\mathbf{n}} = (\hat{\mathbf{n}} b_{\mathbf{n}}^d)^{1/2} \mathbf{c}^\tau (W_{\mathbf{n}}(\mathbf{x}) - \mathrm{E} W_{\mathbf{n}}(\mathbf{x})) = A_{\mathbf{n}} - \mathrm{E} A_{\mathbf{n}}.$$

Now, let us decompose $\hat{\mathbf{n}}^{-1/2} S_{\mathbf{n}}$ into smaller pieces involving "large" and "small" blocks. More specifically, consider [all sums run over $\mathbf{i} := (i_1, \ldots, i_N)$]

$$U(1, \mathbf{n}, \mathbf{x}, \mathbf{j}) := \sum_{\substack{i_k = j_k(p_k+q)+1 \\ k=1,\ldots,N}}^{j_k(p_k+q)+p_k} \zeta_{\mathbf{ni}}(\mathbf{x}),$$

$$U(2, \mathbf{n}, \mathbf{x}, \mathbf{j}) := \sum_{\substack{i_k = j_k(p_k+q)+1 \\ k=1,\ldots,N-1}}^{j_k(p_k+q)+p_k} \sum_{i_N = j_N(p_N+q)+p_N+1}^{(j_N+1)(p_N+q)} \zeta_{\mathbf{ni}}(\mathbf{x}),$$

$$U(3, \mathbf{n}, \mathbf{x}, \mathbf{j}) := \sum_{\substack{i_k = j_k(p_k+q)+1 \\ k=1,\ldots,N-2}}^{j_k(p_k+q)+p_k} \sum_{i_{N-1} = j_{N-1}(p_{N-1}+q)+p_{N-1}+1}^{(j_{N-1}+1)(p_{N-1}+q)} \sum_{i_N = j_N(p_N+q)+1}^{j_N(p_N+q)+p_N} \zeta_{\mathbf{ni}}(\mathbf{x}),$$

$$U(4, \mathbf{n}, \mathbf{x}, \mathbf{j}) := \sum_{\substack{i_k = j_k(p_k+q)+1 \\ k=1,\ldots,N-2}}^{j_k(p_k+q)+p_k} \sum_{i_{N-1} = j_{N-1}(p_{N-1}+q)+p_{N-1}+1}^{(j_{N-1}+1)(p_{N-1}+q)} \sum_{i_N = j_N(p_N+q)+p_N+1}^{(j_N+1)(p_N+q)} \zeta_{\mathbf{ni}}(\mathbf{x}),$$

and so on. Note that

$$U(2^N - 1, \mathbf{n}, \mathbf{x}, \mathbf{j}) := \sum_{\substack{i_k = j_k(p_k+q)+p_k+1 \\ k=1,\ldots,N-1}}^{(j_k+1)(p_k+q)} \sum_{i_N = j_N(p_N+q)+1}^{j_N(p_N+q)+p_N} \zeta_{\mathbf{ni}}(\mathbf{x})$$

and

$$U(2^N, \mathbf{n}, \mathbf{x}, \mathbf{j}) := \sum_{\substack{i_k = j_k(p_k+q)+p_k+1 \\ k=1,\ldots,N}}^{(j_k+1)(p_k+q)} \zeta_{\mathbf{ni}}(\mathbf{x}).$$

Without loss of generality, assume that, for some integers $r_1, \ldots, r_N$, $\mathbf{n} = (n_1, \ldots, n_N)$ is such that $n_1 = r_1(p_1 + q), \ldots, n_N = r_N(p_N + q)$, with $r_k \to \infty$ for all $k = 1, \ldots, N$. For each integer $1 \leq i \leq 2^N$, define

$$T(\mathbf{n}, \mathbf{x}, i) := \sum_{\substack{j_k=0 \\ k=1,\ldots,N}}^{r_k-1} U(i, \mathbf{n}, \mathbf{x}, \mathbf{j}).$$

Clearly, $S_\mathbf{n} = \sum_{i=1}^{2^N} T(\mathbf{n}, \mathbf{x}, i)$. Note that $T(\mathbf{n}, \mathbf{x}, 1)$ is the sum of the random variables $\zeta_{\mathbf{ni}}$ over "large" blocks, whereas $T(\mathbf{n}, \mathbf{x}, i), 2 \leq i \leq 2^N$, are sums over "small" blocks. If it is not the case that $n_1 = r_1(p_1+q), \ldots, n_N = r_N(p_N+q)$ for some integers $r_1, \ldots, r_N$, then an additional term $T(\mathbf{n}, \mathbf{x}, 2^N + 1)$, say, containing all the $\zeta_{\mathbf{nj}}$'s that are not included in the big or small blocks, can be considered. This term will not change the proof much. The general approach consists in showing that, as $\mathbf{n} \to \infty$,

$$(5.18) \quad Q_1 := \left| \mathrm{E}[\exp[iuT(\mathbf{n}, \mathbf{x}, 1)]] - \prod_{\substack{j_k=0 \\ k=1,\ldots,N}}^{r_k-1} \mathrm{E}[\exp[iuU(1, \mathbf{n}, \mathbf{x}, \mathbf{j})]] \right| \to 0,$$

$$(5.19) \quad Q_2 := \hat{\mathbf{n}}^{-1} \mathrm{E}\left( \sum_{i=2}^{2^N} T(\mathbf{n}, \mathbf{x}, i) \right)^2 \to 0,$$

$$(5.20) \quad Q_3 := \hat{\mathbf{n}}^{-1} \sum_{\substack{j_k=0 \\ k=1,\ldots,N}}^{r_k-1} \mathrm{E}[U(1, \mathbf{n}, \mathbf{x}, \mathbf{j})]^2 \to \sigma^2,$$

$$(5.21) \quad Q_4 := \hat{\mathbf{n}}^{-1} \sum_{\substack{j_k=0 \\ k=1,\ldots,N}}^{r_k-1} \mathrm{E}[(U(1, \mathbf{n}, \mathbf{x}, \mathbf{j}))^2 I\{|U(1, \mathbf{n}, \mathbf{x}, \mathbf{j})| > \varepsilon\sigma\hat{\mathbf{n}}^{1/2}\}] \to 0,$$

for every $\varepsilon > 0$. Note that

$$[A_\mathbf{n} - \mathrm{E}A_\mathbf{n}]/\sigma = (\hat{\mathbf{n}}b_\mathbf{n}^d)^{1/2}\mathbf{c}^\tau[W_\mathbf{n}(\mathbf{x}) - \mathrm{E}W_\mathbf{n}(\mathbf{x})]/\sigma = S_\mathbf{n}/(\sigma\hat{\mathbf{n}}^{1/2})$$

$$= T(\mathbf{n}, \mathbf{x}, 1)/(\sigma\hat{\mathbf{n}}^{1/2}) + \sum_{i=2}^{2^N} T(\mathbf{n}, \mathbf{x}, i)/(\sigma\hat{\mathbf{n}}^{1/2}).$$

The term $\sum_{i=2}^{2^N} T(\mathbf{n}, \mathbf{x}, i)/(\sigma\hat{\mathbf{n}}^{1/2})$ is asymptotically negligible by (5.19). The random variables $U(1, \mathbf{n}, \mathbf{x}, \mathbf{j})$ are asymptotically mutually independent by (5.18). The asymptotic normality of $T(\mathbf{n}, \mathbf{x}, 1)/(\sigma \hat{\mathbf{n}}^{1/2})$ follows from (5.20) and the Lindeberg–Feller condition (5.21). The lemma thus follows if we can prove (5.18)–(5.21). This proof is given here. The arguments are reminiscent of those used by Masry (1986) and Nakhapetyan (1987).

Before turning to the end of the proof of Lemma 3.1, we establish the following preliminary lemma, which significantly reinforces Lemma 3.1 in Tran (1990).

LEMMA 5.3. *Let the spatial process $\{Y_\mathbf{i}, \mathbf{X}_\mathbf{i}\}$ satisfy the mixing property (2.1), and denote by $\widetilde{U}_j$, $j = 1, \ldots, M$, an $M$-tuple of measurable functions such that $\widetilde{U}_j$ is measurable with respect to $\{(Y_\mathbf{i}, \mathbf{X}_\mathbf{i}), \mathbf{i} \in \tilde{\mathcal{I}}_j\}$, where $\tilde{\mathcal{I}}_j \subset \mathcal{I}_\mathbf{n}$. If $\mathrm{Card}(\tilde{\mathcal{I}}_j) \leq p$ and $d(\tilde{\mathcal{I}}_\ell, \tilde{\mathcal{I}}_j) \geq q$ for any $\ell \neq j$, then*

$$\left| \mathrm{E}\left[ \exp\left\{ iu \sum_{j=1}^M \widetilde{U}_j \right\} \right] - \prod_{j=1}^M \mathrm{E}[\exp\{iu\widetilde{U}_j\}] \right| \leq C \sum_{j=1}^{M-1} \psi(p, (M-j)p)\varphi(q),$$

*where $i = \sqrt{-1}$.*

PROOF. Let $a_j := \exp\{iu\widetilde{U}_j\}$. Then

$$\mathrm{E}[a_1 \cdots a_M] - \mathrm{E}[a_1] \cdots \mathrm{E}[a_M]$$
$$= \mathrm{E}[a_1 \cdots a_M] - \mathrm{E}[a_1]\mathrm{E}[a_2 \cdots a_M]$$
$$+ \mathrm{E}[a_1]\{\mathrm{E}[a_2 \cdots a_M] - \mathrm{E}[a_2]\mathrm{E}[a_3 \cdots a_M]\}$$
$$+ \cdots + \mathrm{E}[a_1]\mathrm{E}[a_2] \cdots \mathrm{E}[a_{M-2}]\{\mathrm{E}[a_{M-1}a_M] - \mathrm{E}[a_{M-1}]\mathrm{E}[a_M]\}.$$

Since $|\mathrm{E}[a_i]| \leq 1$,

$$|\mathrm{E}[a_1 \cdots a_M] - \mathrm{E}[a_1] \cdots \mathrm{E}[a_M]|$$
$$\leq |\mathrm{E}[a_1 \cdots a_M] - \mathrm{E}[a_1]\mathrm{E}[a_2 \cdots a_M]|$$
$$+ |\mathrm{E}[a_2 \cdots a_M] - \mathrm{E}[a_2]\mathrm{E}[a_3 \cdots a_M]|$$
$$+ \cdots + |\mathrm{E}[a_{M-1}a_M] - \mathrm{E}[a_{M-1}]\mathrm{E}[a_M]|.$$

Note that $d(I_\ell, I_j) \geq q$ for any $\ell \neq j$. The lemma then follows by applying Lemma 5.1(ii) to each term on the right-hand side. □

PROOF OF LEMMA 3.1 (continued). In order to complete the proof of Lemma 3.1, we still have to prove (5.18)–(5.21).

PROOF OF (5.18). Ranking the random variables $U(1, \mathbf{n}, \mathbf{x}, \mathbf{j})$ in an arbitrary manner, refer to them as $\widetilde{U}_1, \ldots, \widetilde{U}_M$. Note that $M = \prod_{k=1}^N r_k = \hat{\mathbf{n}}\{\prod_{k=1}^N (p_k + q)\}^{-1} \leq \hat{\mathbf{n}}/p$, where $p = \prod_{k=1}^N p_k$. Let

$$\mathcal{I}(1, \mathbf{n}, \mathbf{x}, \mathbf{j}) := \{\mathbf{i} : j_k(p_k + q) + 1 \leq i_k \leq j_k(p_k + q) + p_k, k = 1, \ldots, N\}.$$

The distance between two distinct sets $\mathcal{I}(1,\mathbf{n},\mathbf{x},\mathbf{j})$ and $\mathcal{I}(1,\mathbf{n},\mathbf{x},\mathbf{j}')$ is at least $q$. Clearly, $\mathcal{I}(1,\mathbf{n},\mathbf{x},\mathbf{j})$ is the set of sites involved in $U(1,\mathbf{n},\mathbf{x},\mathbf{j})$. As for the set of sites $\widetilde{\mathcal{I}}_j$ associated with $\widetilde{U}_j$, it contains $p$ elements. Hence, in view of Lemma 5.3 and assumption (A4'),

$$Q_1 \leq C \sum_{k=1}^{M-1} \min\{p,(M-k)p\}\varphi(q) \leq CMp\varphi(q) \leq C\hat{\mathbf{n}}\varphi(q),$$

which tends to zero by condition (B2). $\square$

PROOF OF (5.19). In order to prove (5.19), it is enough to show that

$$\hat{\mathbf{n}}^{-1}\mathrm{E}[T^2(\mathbf{n},\mathbf{x},i)] \to 0 \qquad \text{for any } 2 \leq i \leq 2^N.$$

Without loss of generality, consider $\mathrm{E}[T^2(\mathbf{n},\mathbf{x},2)]$. Ranking the random variables $U(2,\mathbf{n},\mathbf{x},\mathbf{j})$ in an arbitrary manner, refer to them as $\widehat{U}_1,\ldots,\widehat{U}_M$. We have

$$(5.22) \qquad \mathrm{E}[T^2(\mathbf{n},\mathbf{x},2)] = \sum_{i=1}^{M} \mathrm{Var}(\widehat{U}_i) + 2 \sum_{1 \leq i < j \leq M} \mathrm{Cov}(\widehat{U}_i,\widehat{U}_j)$$
$$:= \widehat{V}_1 + \widehat{V}_2 \qquad \text{say.}$$

Since $\mathbf{X}_n$ is stationary [recall that $\zeta_{\mathbf{nj}}(\mathbf{x}) := b_{\mathbf{n}}^{-d/2}\Delta_{\mathbf{j}}(\mathbf{x})$],

$$\mathrm{Var}(\widehat{U}_i) = \mathrm{E}\!\left[\left(\sum_{\substack{i_k=1\\k=1,\ldots,N-1}}^{p_k} \sum_{i_N=1}^{q} \zeta_{\mathbf{ni}}(\mathbf{x})\right)^{\!2}\right] + \sum_{\mathbf{i} \neq \mathbf{j} \in \mathcal{J}} \mathrm{E}[\zeta_{\mathbf{nj}}(\mathbf{x})\zeta_{\mathbf{ni}}(\mathbf{x})] := \widehat{V}_{11} + \widehat{V}_{12},$$

where $\mathcal{J} = \mathcal{J}(\mathbf{p},q) := \{\mathbf{i},\mathbf{j} : 1 \leq i_k, j_k \leq p_k, k = 1,\ldots,N-1, \text{ and } 1 \leq i_N, j_N \leq q\}$. From (5.13) and the Lebesgue density theorem [see Chapter 2 of Devroye and Györfi (1985)],

$$\widehat{V}_{11} = \left(\prod_{k=1}^{N-1} p_k\right) q \, \mathrm{Var}\{\zeta_{\mathbf{ni}}(\mathbf{x})\} = \left(\prod_{k=1}^{N-1} p_k\right) q\{b_{\mathbf{n}}^{-d}\mathrm{E}\Delta_{\mathbf{i}}^2(\mathbf{x})\} \leq C\left(\prod_{k=1}^{N-1} p_k\right) q.$$

Thus, applying Lemma 5.2 with $n_k = p_k$, $k = 1,\ldots,N-1$, and $n_N = q$ yields

$$\widehat{V}_{12} = b_{\mathbf{n}}^{-d} \sum_{\mathbf{i} \neq \mathbf{j} \in \mathcal{J}} \mathrm{E}[\Delta_{\mathbf{j}}(\mathbf{x})\Delta_{\mathbf{i}}(\mathbf{x})]$$
$$\leq Cb_{\mathbf{n}}^{-d}\left[b_{\mathbf{n}}^{\delta d/(4+\delta)+d}\left(\prod_{k=1}^{N-1} p_k c_{\mathbf{n}k}\right) qc_{\mathbf{n}N}\right.$$
$$\left.+ b_{\mathbf{n}}^{2d/(2+\delta)}\left(\prod_{k=1}^{N-1} p_k\right) q \sum_{k=1}^{N} \sum_{t=c_{\mathbf{n}k}}^{\|\mathbf{n}\|} t^{N-1}\{\varphi(t)\}^{\delta/(2+\delta)}\right]$$

$$= C\left(\prod_{k=1}^{N-1} p_k\right) q \left[ b_{\mathbf{n}}^{\delta d/(4+\delta)} \left(\prod_{k=1}^{N} c_{\mathbf{n}k}\right) \right.$$
$$\left. + b_{\mathbf{n}}^{-\delta d/(2+\delta)} \sum_{k=1}^{N} \sum_{t=c_{\mathbf{n}k}}^{\infty} t^{N-1} \{\varphi(t)\}^{\delta/(2+\delta)} \right]$$
$$:= C\left(\prod_{k=1}^{N-1} p_k\right) q \pi_{\mathbf{n}}.$$

It follows that

$$\hat{\mathbf{n}}^{-1} \widehat{V}_1 = \hat{\mathbf{n}}^{-1} M(\hat{V}_{11} + \hat{V}_{12})$$
(5.23)
$$\leq \hat{\mathbf{n}}^{-1} MC\left(\prod_{k=1}^{N-1} p_k\right) q[1 + \pi_{\mathbf{n}}] \leq C(q/p_N)[1 + \pi_{\mathbf{n}}].$$

Set

$$\mathcal{I}(2, n, \mathbf{x}, \mathbf{j}) := \{\mathbf{i} : j_k(p_k + q) + 1 \leq i_k \leq j_k(p_k + q) + p_k, 1 \leq k \leq N-1,$$
$$j_N(p_N + q) + p_N + 1 \leq i_N \leq (j_N + 1)(p_N + q)\}.$$

Then $U(2, \mathbf{n}, \mathbf{x}, \mathbf{j}) = \sum_{\mathbf{i} \in \mathcal{I}(2,\mathbf{n},\mathbf{x},\mathbf{j})} \zeta_{\mathbf{n}\mathbf{i}}(\mathbf{x})$. Since $p_k > q$, if $\mathbf{i}$ and $\mathbf{i}'$ belong to two distinct sets $\mathcal{I}(2, \mathbf{n}, \mathbf{x}, \mathbf{j})$ and $\mathcal{I}(2, \mathbf{n}, \mathbf{x}, \mathbf{j}')$, then $\|\mathbf{i} - \mathbf{i}'\| > q$. In view of (5.8) and (5.22), we obtain

$$|\widehat{V}_2| \leq C \sum_{\{\mathbf{i},\mathbf{j}:\|\mathbf{i}-\mathbf{j}\| \geq q,\, 1 \leq i_k, j_k \leq n_k\}} |\mathrm{E}[\zeta_{\mathbf{n}\mathbf{i}}(\mathbf{x}) \zeta_{\mathbf{n}\mathbf{j}}(\mathbf{x})]|$$
$$\leq C b_{\mathbf{n}}^{-d} \sum_{\{\mathbf{i},\mathbf{j}:\|\mathbf{i}-\mathbf{j}\| \geq q,\, 1 \leq i_k, j_k \leq n_k\}} |\mathrm{E}[\Delta_{\mathbf{n}\mathbf{i}}(\mathbf{x}) \Delta_{\mathbf{n}\mathbf{j}}(\mathbf{x})]|$$
$$\leq C b_{\mathbf{n}}^{-d} \sum_{\{\mathbf{i},\mathbf{j}:\|\mathbf{i}-\mathbf{j}\| \geq q,\, 1 \leq i_k, j_k \leq n_k\}} b_{\mathbf{n}}^{2d/(2+\delta)} \{\varphi(\|\mathbf{j} - \mathbf{i}\|)\}^{\delta/(2+\delta)}$$
(5.24)
$$\leq C b_{\mathbf{n}}^{-\delta d/(2+\delta)} \left(\prod_{k=1}^{N} n_k\right) \left(\sum_{t=q}^{\|\mathbf{n}\|} t^{N-1} \{\varphi(t)\}^{\delta/(2+\delta)}\right).$$

Take $c_{\mathbf{n}k}^a = b_{\mathbf{n}}^{-\delta d/(2+\delta)} \to \infty$. Condition (B3) implies that $q b_{\mathbf{n}}^{\delta d/a(2+\delta)} > 1$, so that $c_{\mathbf{n}k} < q \leq p_k$. Then, as proved in (5.15) and (5.16), it follows from assumption (A4) that $\pi_{\mathbf{n}} \to 0$. Thus, from (5.22), (5.23) and (5.24),

$$\hat{\mathbf{n}}^{-1} \mathrm{E}[T^2(\mathbf{n}, \mathbf{x}, 2)] \leq C(q/p_N)[1 + \pi_{\mathbf{n}}] + C b_{\mathbf{n}}^{-\delta d/(2+\delta)} \left(\sum_{t=q}^{\infty} t^{N-1} \{\varphi(t)\}^{\delta/(2+\delta)}\right),$$

which tends to zero by $q/p_N \to 0$ and condition (B3); (5.19) follows. □

PROOF OF (5.20). Let $S'_\mathbf{n} := T(\mathbf{n}, \mathbf{x}, 1)$ and $S''_\mathbf{n} := \sum_{i=2}^{2^N} T(\mathbf{n}, \mathbf{x}, i)$. Then $S'_\mathbf{n}$ is a sum of $Y_\mathbf{j}$'s over the "large" blocks, $S''_\mathbf{n}$ over the "small" ones. Lemma 3.2 implies $\hat{\mathbf{n}}^{-1} \mathrm{E}[|S_\mathbf{n}|^2] \to \sigma^2$. This, combined with (5.19), entails $\hat{\mathbf{n}}^{-1} \mathrm{E}[|S'_\mathbf{n}|^2] \to \sigma^2$. Now,

$$
(5.25) \quad \hat{\mathbf{n}}^{-1} \mathrm{E}[|S'_\mathbf{n}|^2] = \hat{\mathbf{n}}^{-1} \sum_{\substack{j_k=0 \\ k=1,\ldots,N}}^{r_k-1} \mathrm{E}[U^2(1, \mathbf{n}, \mathbf{x}, \mathbf{j})] + \hat{\mathbf{n}}^{-1} \sum_{\mathbf{i} \neq \mathbf{j} \in \mathcal{J}^*} \mathrm{Cov}(U(1, \mathbf{n}, \mathbf{x}, \mathbf{j}), U(1, \mathbf{n}, \mathbf{x}, \mathbf{i})),
$$

where $\mathcal{J}^* = \mathcal{J}^*(\mathbf{p}, q) := \{\mathbf{i}, \mathbf{j} : 1 \leq i_k, j_k \leq r_k - 1, k = 1, \ldots, N\}$. Observe that (5.20) follows from (5.25) if the last sum in the right-hand side of (5.25) tends to zero as $\mathbf{n} \to \infty$. Using the same argument as in the derivation of the bound (5.22) for $\widehat{V}_2$, this sum can be bounded by

$$
Cb_\mathbf{n}^{-\delta d/(2+\delta)} \sum_{\|\mathbf{i}\|>q} \sum_{\substack{i_k=1 \\ k=1,\ldots,N}}^{n_k-1} \{\varphi(\|\mathbf{i}\|)\}^{\delta/(2+\delta)} \leq Cb_\mathbf{n}^{-\delta d/(2+\delta)} \left( \sum_{t=q}^{\infty} t^{N-1} \{\varphi(t)\}^{\delta/(2+\delta)} \right),
$$

which tends to zero by condition (B3). □

PROOF OF (5.21). We need a truncation argument because $Z_\mathbf{i}$ is not necessarily bounded. Set $Z_\mathbf{i}^L := Z_\mathbf{i} I_{\{|Z_\mathbf{i}| \leq L\}}$, $\eta_\mathbf{i}^L := Z_\mathbf{i}^L K_\mathbf{c}((\mathbf{X}_\mathbf{i} - \mathbf{x})/b_\mathbf{n})$, $\Delta_\mathbf{i}^L := \eta_\mathbf{i}^L - \mathrm{E}\eta_\mathbf{i}^L$, $\zeta_{\mathbf{ni}}^L := b_\mathbf{n}^{-d/2} \Delta_\mathbf{i}^L$, where $L$ is a fixed positive constant, and define $U^L(1, \mathbf{n}, \mathbf{x}, \mathbf{j}) := \sum_{\mathbf{i} \in \mathcal{I}(1,\mathbf{n},\mathbf{x},\mathbf{j})} \zeta_{\mathbf{ni}}^L$. Put

$$
Q_4^L := \hat{\mathbf{n}}^{-1} \sum_{\substack{j_k=0 \\ k=1,\ldots,N}}^{r_k-1} \mathrm{E}[(U^L(1, \mathbf{n}, \mathbf{x}, \mathbf{j}))^2 I\{|U^L(1, \mathbf{n}, \mathbf{x}, \mathbf{j})| > \varepsilon \sigma \hat{\mathbf{n}}^{1/2}\}].
$$

Clearly, $|\zeta_{\mathbf{ni}}^L| \leq CLb_\mathbf{n}^{-d/2}$. Therefore $|U^L(1, \mathbf{n}, \mathbf{x}, \mathbf{j})| < CLpb_\mathbf{n}^{-d/2}$. Hence

$$
Q_4^L \leq C\hat{\mathbf{p}}^2 b_\mathbf{n}^{-d} \hat{\mathbf{n}}^{-1} \sum_{\substack{j_k=0 \\ k=1,\ldots,N}}^{r_k-1} \mathrm{P}[U^L(1, \mathbf{n}, \mathbf{x}, \mathbf{j}) > \varepsilon \sigma \hat{\mathbf{n}}^{1/2}].
$$

Now, $U^L(1, \mathbf{n}, \mathbf{x}, \mathbf{j})/(\sigma \hat{\mathbf{n}}^{1/2}) \leq C\hat{\mathbf{p}}(\hat{\mathbf{n}} b_\mathbf{n}^d)^{-1/2} \to 0$, since $\hat{\mathbf{p}} = [(\hat{\mathbf{n}} b_\mathbf{n}^d)^{1/2}/s_\mathbf{n}]$, where $s_\mathbf{n} \to \infty$. Thus $\mathrm{P}[U^L(1, \mathbf{n}, \mathbf{x}, \mathbf{j}) > \varepsilon \sigma \hat{\mathbf{n}}^{1/2}] = 0$ at all $\mathbf{j}$ for sufficiently large $\hat{\mathbf{n}}$. Thus $Q_4^L = 0$ for large $\hat{\mathbf{n}}$, and (5.21) holds for the truncated variables. Hence

$$
(5.26) \quad \hat{\mathbf{n}}^{-1/2} S_\mathbf{n}^L := \hat{\mathbf{n}}^{-1/2} \sum_{\substack{j_k=1 \\ k=1,\ldots,N}}^{n_k} \zeta_{\mathbf{nj}}^L \xrightarrow{\mathcal{L}} N(0, \sigma_L^2),
$$

30    M. HALLIN, Z. LU AND L. T. TRAN
where $\sigma_L^2 := \mathrm{Var}(Z_\mathbf{i}^L|\mathbf{X_i} = \mathbf{x})f(\mathbf{x}) \int K_\mathbf{c}^2(\mathbf{u})\,d\mathbf{u}$.

Defining $S_\mathbf{n}^{L*} := \sum_{j_k=1; k=1,\ldots,N}^{n_k}(\zeta_\mathbf{nj} - \zeta_\mathbf{nj}^L)$, we have $S_\mathbf{n} = S_\mathbf{n}^L + S_\mathbf{n}^{L*}$. Note that

$$\begin{aligned}
|\mathrm{E}[\exp(iuS_\mathbf{n}/\hat{\mathbf{n}}^{1/2})] - \exp(-u^2\sigma^2/2)| \\
\leq |\mathrm{E}[\exp(iuS_\mathbf{n}^L/\hat{\mathbf{n}}^{1/2}) - \exp(-u^2\sigma_L^2/2)]\exp(iuS_\mathbf{n}^{L*}/\hat{\mathbf{n}}^{1/2})| \\
+ |\mathrm{E}[\exp(iuS_\mathbf{n}^{L*}/\hat{\mathbf{n}}^{1/2}) - 1]\exp(-u^2\sigma_L^2/2)| \\
+ |\exp(-u^2\sigma_L^2/2) - \exp(-u^2\sigma^2/2)| \\
= E_1 + E_2 + E_3, \qquad \text{say.}
\end{aligned}$$

Letting $\mathbf{n} \to \infty$, $E_1$ tends to zero by (5.26) and the dominated convergence theorem. Letting $L$ go to infinity, the dominated convergence theorem also implies that $\sigma_L^2 := \mathrm{Var}(Z_\mathbf{i}^L|\mathbf{X_i} = \mathbf{x})f(\mathbf{x}) \int K_\mathbf{c}^2(\mathbf{u})\,d\mathbf{u}$ converges to

$$\mathrm{Var}(Z_\mathbf{i}|\mathbf{X_i} = \mathbf{x})f(\mathbf{x})\int K_\mathbf{c}^2(\mathbf{u})\,d\mathbf{u} = \mathrm{Var}(Y_\mathbf{i}|\mathbf{X_i} = \mathbf{x})f(\mathbf{x})\int K_\mathbf{c}^2(\mathbf{u})\,d\mathbf{u} := \sigma^2,$$

and hence that $E_3$ tends to zero. Finally, in order to prove that $E_2$ also tends to zero, it suffices to show that $S_\mathbf{n}^{L*}/\hat{\mathbf{n}}^{1/2} \to 0$ in probability as first $\mathbf{n} \to \infty$ and then $L \to \infty$, which in turn would follow if we could show that

$$\mathrm{E}[(S_\mathbf{n}^{L*}/\hat{\mathbf{n}}^{1/2})^2] \to \mathrm{Var}(|Z_\mathbf{i}|I_{\{|Z_\mathbf{i}|>L\}}|\mathbf{X_i} = \mathbf{x})f(\mathbf{x})\int K_\mathbf{c}^2(\mathbf{u})\,d\mathbf{u} \qquad \text{as } \mathbf{n} \to \infty.$$

This follows along the same lines as Lemma 3.2. □

The proof of Lemma 3.1 is thus complete. □

**Acknowledgments.** The authors are grateful to two anonymous referees and an Associate Editor for their careful reading of the first version of this paper, and for their insightful and constructive comments.
## REFERENCES

Anselin, L. and Florax, R. J. G. M. (1995). *New Directions in Spatial Econometrics*. Springer, Berlin.

Besag, J. E. (1974). Spatial interaction and the statistical analysis of lattice systems (with discussion). *J. Roy. Statist. Soc. Ser. B* **36** 192–236. MR373208

Biau, G. (2003). Spatial kernel density estimation. *Math. Methods Statist.* **12** 371–390. MR2054154

Biau, G. and Cadre, B. (2004). Nonparametric spatial prediction. *Stat. Inference Stoch. Process.* **7** 327–349. MR2111294

Boente, G. and Fraiman, R. (1988). Consistency of a nonparametric estimate of a density function for dependent variables. *J. Multivariate Anal.* **25** 90–99. MR935296

Bolthausen, E. (1982). On the central limit theorem for stationary mixing random fields. *Ann. Probab.* **10** 1047–1050. MR672305

M. Hallin  
Institut de Statistique et  
  de Recherche Opérationnelle  
Université Libre de Bruxelles  
Campus de la Plaine CP 210  
B-1050 Bruxelles  
Belgium  
e-mail: mhallin@ulb.ac.be

Z. Lu  
Institute of Systems Science  
Academy of Mathematics  
  and Systems Sciences  
Chinese Academy of Sciences  
Beijing  
China  
and  
Department of Statistics  
London School of Economics  
Houghton Street  
London WC2A 2AE  
United Kingdom  
e-mail: z.lu@lse.ac.uk

L. T. Tran  
Department of Mathematics  
Rawles Hall  
Indiana University  
Bloomington, Indiana 47405  
USA  
e-mail: tran@indiana.edu